\documentclass[11pt]{amsart}
\usepackage{amssymb, amsmath, amsthm}
\usepackage[colorlinks=true,linkcolor=blue,citecolor=red]{hyperref}
\usepackage{color}
\usepackage{epsfig}
\usepackage{amsmath}
\usepackage{amssymb}
\usepackage{amscd}
\usepackage{graphicx}
\usepackage{mathrsfs}
\usepackage[T1]{fontenc}
\numberwithin{equation}{section}
\newtheorem{lemma}{Lemma}[section]
   \newtheorem{theorem}{Theorem}
   \newtheorem{re}{Remark}[section]
   \newtheorem{prop}{Proposition}[section]

\def\spn{\mathop{\mathrm{span}}}

\usepackage{epsfig}

\numberwithin{equation}{section}
\numberwithin{theorem}{section}
\numberwithin{prop}{section}
\numberwithin{lemma}{section}
\numberwithin{re}{section}
\numberwithin{coro}{section}

\DeclareMathOperator{\sech}{sech}

\newcommand{\al}{\alpha}
\newcommand{\ppt}{\phi_{c_j(t)}}
\newcommand{\ppo}{\phi_{c_j(0)}}

\newcommand{\e}{\epsilon}
\newcommand{\ka}{\kappa}
\newcommand{\R}{\mathbb{R}}

\newcommand{\Z}{\mathbb{Z}}

\newcommand{\cau}{\mathcal{U}}
\newcommand{\cav}{\mathcal{V}}

\newcommand{\cas}{\mathcal{S}}
\newcommand{\cal}{\mathcal{L}}

\newcommand{\aaa}{\int_{\R}}
\newcommand{\rmd}{\textrm{d}}



\thanks{ Email: \dag liji@hust.edu.cn, \ddag yliu@uta.edu, \S {wuq@ohio.edu}. }

\begin{document}

\title[Orbital stability of sum of smooth DP solitons]{Orbital Stability of the sum of Smooth solitons in the Degasperis-Procesi Equation}


\author[Li]{Ji Li}
\address{Ji Li
\newline
School of Mathematics and Statistics, Huazhong University of Science and Technology, Wuhan 430074, P. R. China}
\email{liji@hust.edu.cn}

\author[Liu]{Yue Liu}
\address{Yue Liu \newline
Department of Mathematics, University of Texas at Arlington, TX 76019}
\email{yliu@uta.edu}

\author[Wu]{Qiliang Wu}

\address{Qiliang Wu (Corresponding author)\newline
Department of Mathematics, Ohio University, Athen, OH 45701}
\email{wuq@ohio.edu}

\begin{abstract}
{The Degasperis-Procesi (DP) equation is an integrable Camassa-Holm-type model as an asymptotic approximation for the unidirectional propagation of shallow water waves. 
This work is to establish the $L^2\cap L^\infty$ orbital stability  of  a wave train containing $N$ smooth solitons which are well separated. The main difficulties stem from the subtle nonlocal structure of the DP equation. One consequence is that the energy space of the DE equation based on the conserved quantity induced by the translation symmetry is only equivalent to the $L^2$-norm, which by itself can not bound the higher-order nonlinear terms in the Lagrangian. Our remedy is to introduce \textit{a priori } estimates based on certain smooth initial conditions. Moreover, another consequence is that the nonlocal structure of the DP equation significantly complicates the verification of the monotonicity of local momentum and the positive definiteness of a refined quadratic form of the orthogonalized perturbation.
}
\end{abstract}

\maketitle

\section{Introduction}
{Solitary wave phenomena were first scientifically described and studied by Russell \cite{Rus} in 1834. One significant observation of Russell, among others, was the existence of solitary waves whose shapes do not disperse as they propagate. The mathematical vindication of this controversial result was pursued by many mathematicians such as  Boussinesq \cite{Bou}, Rayleigh \cite{Ray} and finally resolved in 1895 by Korteweg and de Vries \cite{KdV} whom derived a nonlinear PDE, now called the Korteweg-de Vries (KdV) equation, taking its dimensionless form,
\begin{equation}
u_t+6uu_x+u_{xxx}=0,
\end{equation}
and showed that the KdV equation admits solitary wave solution
\[
\phi(\xi;c)=\frac{c}{2}\sech^2\left( \frac{\sqrt{c}}{2}\xi \right), \quad \text{ where }\xi\triangleq x-ct.
\]
Except for that, such solitary wave phenomena were largely dismissed till the 1960s when the remarkable property, that these solitary waves can preserve their shape and velocity after crossing each other, thus named as solitons, was numerically observed by Zabusky and Kruskal \cite{ZK}, and later theoretically explained by Gardner \textit{et. al.} \cite{GGKM} via the inverse scattering method. This inverse scattering method announced the advent of modern study of solitary waves, which eventually unveiled the universality and significance of solitons, and established a research field vibrant till this day; see \cite{Miu, New} for more detailed history and general exposition. One essential result from the inverse scattering method is the existence of $N$-solitons, which are explicit solutions in the form of nonlinear superpositions of $N$ single solitons, asymptotically reducing to linear superposition of $N$ single solitons as $t\to\infty$. The inverse scattering method has since been adapted and generalized to solve many other completely integrable infinite dimensional systems, typically accommodating $N$-solitons. Various other techniques, including the Hirota's direct method and the Wronskian method, are also developed for the derivation of $N$-solitons; see \cite{Hiro, HJ, NF} and the references therein for more details.}

{It is also vital to check the stability of $N$-solitons for the validity of mathematical models of these physical systems bearing soliton phenomena. While for single solitary waves ($N=1$) case, there are classical Grillakis-Shatah-Strauss theory, Zakharov-Shabat's inverse scattering method and so on, the development in stability theory for $N$-solitons is less systematic. Similar to the $n=1$ case, Maddocks and Sachs \cite{MS} showed that $N$ solitons of the KdV equation are non-isolated constrained minimizers and thus stable in $H^N(\R)$.  This result was improved by Martel, Merle and Tsai whom obtained stability of the sum of $N$ single solitons, instead of $N$-solitons, in the energy space of subcritical generalized KdV equation \cite{MMT} and of nonlinear Schr\"{o}dinger equations \cite{MMT-2}, and this is the path we follow in this paper. Their method was later extended to many other PDEs such as the Gross-Pitaevskii equation \cite{BGS}, the Landau-Lifshitz equation \cite{deL-G} and the nonlinear Schr\"{o}dinger equation with derivative cubic nonlinearity \cite{SW}.  In addition, stability of $N$-solitons was also proved in weighted spaces \cite{Per, RSS}. For instability of $N$-solitons results we refer to \cite{DM, BSS} and references therein.}

{In this paper, we study the stability of the sum of $N$ smooth solitary waves for the Degasperis-Procesi (DP) Equation
\begin{equation}\label{DPk}
m_t+2\kappa u_x+3mu_x+um_x=0,\quad x\in\R,  \;  t>0,
\end{equation}
where $m\triangleq u-u_{xx}$  is the momentum density and $ \kappa> 0 $ is a parameter related to the critical shallow water speed. It is noted  that the DP equation, the Camassa-Holm (CH) equation \cite{CH, F-F}
\begin{equation} \label{CH}
m_t+2\kappa u_x+2mu_x+um_x=0,
\end{equation}
and the KdV equation are the only three integrable candidates in a broad family of third-order dispersive PDEs; see \cite{LLW1} for a detailed comparison of their similarities and differences. We have proved the existence and spectrum stability of smooth solitary waves of the DP equation in \cite{LLW1}. More specifically, }
the DP equation can be written as an infinite dimensional Hamiltonian PDE; that is,
\begin{equation}\label{DP}
u_t=J\frac{\delta H}{\delta u}(u),
\end{equation}
where $$J\triangleq\partial_x(4-\partial_x^2)(1-\partial_x^2)^{-1},\quad H(u)\triangleq-\frac{1}{6}\int \left ( u^3+6\kappa  \left ((4-\partial_x^2)^{-\frac{1}{2}}u \right )^2 \right ) \, dx. $$
Or more explicitly:
\begin{equation}\label{kDP}
\partial_tu+\partial_x\left(\frac{1}{2}u^2+p*(\frac{3}{2}u^2+2\kappa u)\right)=0,\quad t>0, \; \, x\in\R,
\end{equation}
where $p(x)=\frac{1}{2}e^{-\vert x\vert}$ is the impulse response corresponding to the operator $1-\partial_x^2$ so that for all $f\in L^2(\R)$,
\[
(1-\partial_x^2)^{-1}f=p*f.
\]
The translation invariance gives rise to a conserved quantity:
\begin{equation}\label{E3DP}
S(u)=\frac{1}{2}\int_\R u\cdot(1-\partial_x^2)(4-\partial_x^2)^{-1}u\,dx=\frac{1}{2}\int_\R (4\widehat{u}^2+5\widehat{u}_x^2+\widehat{u}_{xx}^2)\,dx,
\end{equation}
where all through the paper, we 
{introduce the notation
\[
\widehat{f}:=(4-\partial_x^2)^{-1}f
\]}
for convenience. Although the DP equation admits infinitely many conserved quantities, this conserved quantity $S$ and the Hamiltonian $H$ are the essential ones relevant to our study of stability. Compared with its counterpart in the CH equation,
\begin{equation}\label{E3CH}
\widetilde{S}(u)=\int_\R \left ( u^2+u_x^2 \right ) \,dx,
\end{equation}
which is equivalent to the square of the $H^1$-norm, the conserved quantity $S$ is only equivalent to the square of the $L^2$-norm, which by itself can not bound the higher-order (the cubic) nonlinear terms in the Lagrangian, leading to the most significant difficulties one has to overcome in order to prove orbital stability result.

We first recall the existence result of smooth solitary waves established in \cite{LLW1}.
\begin{prop}[existence \cite{LLW1}]\label{profile}
Given the physical condition $c>2\kappa >0$, there exists a unique $c-$speed smooth solitary wave, denoted as $\phi_c(x-ct)$, with its profile $\phi_c$ satisfying
\begin{equation}\label{firstintegral}
\frac{1}{2}(c-\phi_c )^2\phi_{c,x}^2=\phi_c^2 (\frac{1}{2}\phi_c^2 -c\phi_c +\frac{2}{3}\kappa\phi_c+\frac{1}{2}c^2-\kappa c).
\end{equation}
$\phi_c(\cdot)$ is even and strictly decreasing on $(0, \infty)$ with an exponential decay rate $e^{-\sqrt{1-\frac{2\kappa}{c}}}$, i.e.
\begin{equation}\label{exp0}
\vert\phi_c(\cdot)\vert\leq Ce^{-\sqrt{1-\frac{2\kappa}{c}}\vert\cdot\vert} .
\end{equation}
\end{prop}

The current work establishes the following stability result of a train of well separated smooth solitons.

\begin{theorem}[{Orbital stability of $N$-train}]\label{Nstability}
Let $0<2\ka<c_1^0<c_2^0<\cdots<c_N^0$. There exist $\gamma_0, \ A_0,\ L_0,\ \alpha_0>0$ such that the following is true: Let $u_0\in H^s(\R)$, $s>3/2$, and assume that there exist $L>L_0, \ \al<\al_0$, and $x_1^0<x_2^0<\cdots<x_N^0$, such that
\begin{equation}\label{initialerror}
\left\Vert u_0-\sum_{j=1}^N\phi_{c_j^0}(\cdot-x_j^0) \right\Vert_{L^2(\R)}\leq\al,\quad
\mbox{ and for }j=1,\cdots,N-1,
\quad x_{j+1}^0>x_j^0+L.
\end{equation}
In addition, $$w_0\triangleq u_0-(u_0)_{xx}+\frac{2\kappa }{3}>0,$$ is a positive Radon measure in the sense that the mapping $f\mapsto \int_\R fw\textrm{d}x$ gives a continuous linear functional on the space of compact-supported continuous scalar functions equipped with the canonical limit topology. Let $u(t)\in C([0,+\infty),\  H^s(\R))$ be the solution of \eqref{CH}. Then, there exist $x_j(t)$ such that
\begin{equation}\label{error}
\sup_{t\in[0,\infty)}\left\Vert u(t,\cdot)-\sum_{j=1}^N\phi_{c_j^0}(\cdot-x_j(t))\right\Vert_{L^2(\R)}\leq A_0\left(\al+e^{-\gamma_0L/2}\right).
\end{equation}
\end{theorem}

{
\begin{re}
The existence of $N$-solutions of the DP equation is given by Matsuno \cite{Mat-1, Mat-2}.
\end{re}}

\begin{re}
Besides Proposition \ref{prolinfty} on global existence of strong solution, there is also a global existence of weak solutions in $L^2$-space given in \cite{E-L-Y}. The regularity requirement in Theorem \ref{Nstability} can be relaxed to
\[
u_0\in L^2(\R), \quad w=u_0-(u_0)_{xx}+\frac{2\kappa }{3} \text{ is a positive Radon measure.}
\]
The peakon case when $\kappa=0$ can be seen in \cite{KM}.
\end{re}

The stability proof follows the framework seminally developed by Martel, et.al. \cite{MMT}, which localizes the well known Lyapunov functional method developed in \cite{GSS} by exploiting monotonicity of local mass. More specifically, there are three major steps. Firstly, the positive definiteness of a constrained quadratic form should be established, typically based on spectral properties of the second variational derivative of a Lagrangian at solitary waves. Secondly, high order nonlinear terms, if there is any, has to be controlled in order to construct the Lyapunov functional (for the single solitary wave), or equivalently, to apply a bootstrap argument. This is done in this work using $L^2-L^\infty$ \textit{a priori} estimates to bound the cubic nonlinear term, inspired by the recent work of Khorbatly and Molinet \cite{KM}. At last, the proof is concluded with a localization of the above-mentioned Lyapunov function via a monotonicity argument of the localized momentum as well as the positive definiteness of a refined quadratic form of the orthogonalized perturbation. 


The remainder of the paper is organized as follows. In section 2 we prove an \textit{a priori} estimate and prove some useful estimates related to DP. In Section 3, we prove Theorem \ref{Nstability} assuming several lemmas. In Section 4, we prove these lemmas.

\bigskip

\noindent{\bf Notations} 
{We collect here important notations adopted throughout the paper. $f$ and $g$ will be used for general auxiliary functions, which may have different meaning in different places. We use $C$ to universally denote generic constants which are independent of $L$ and $\al$. 
We also introduce the notation
the $L^2$-equivalent inner product
$$(u,v)_{\cas}\triangleq (u,(1-\partial_x^2)(4-\partial_x^2)^{-1}v)=(\frac{\delta S}{\delta u}(u), v).$$
}
Integration $\int$ without index always means $\int_{\R} $, and $\Vert\cdot\Vert_{L^p}$ always means $\Vert\cdot\Vert_{L^p(\R)}$ with $p=2,3,\infty$. $\sum$ without index always means $\sum_{j=1}^N$.

\bigskip

\section{Well-posedness and \emph{a priori} estimates}\label{s:2}
\subsection{Well-posedness}
Generically, for evolutionary equations, the well-posedness of initial value problems serves as the precondition for any qualitative study of their dynamics. For clarity, we list relavant well-posedness results of the DP equation with $\kappa>0$, including uniqueness and existence results, 
whose proofs are based upon the ones for the vanishing linear dispersion case, $(\kappa=0)$, without or with mild modifications.

%

A local well-posedness result for the Cauchy problem \eqref{e:DP-IVP} with $\kappa=0$ is obtained in \cite{Y1} via applying Kato's theorem \cite{K}. With exactly the same argument (thus omitted), we have the following local well-posedness result for the Cauchy problem \eqref{e:DP-IVP} with $\kappa>0$.

\begin{prop}[Uniqueness and local existence of strong solutions]\label{p:loc}
Given the initial profile $u_0\in H^s(\R)$ with $s>\frac{3}{2},$ there exist a maximal  time $T=T(u_0)\in(0,\infty]$, independent of the choice of $s$, and a unique solution $u$ to the Cauchy problem
Consider the following Cauchy problem
\begin{equation}\label{e:DP-IVP}
\begin{cases}
\partial_tu+\partial_x\left(\frac{1}{2}u^2+p*(\frac{3}{2}u^2+2\kappa u)\right)=0,\\
u(0,x)=u_0(x).
\end{cases}
\end{equation}
such that \ $u=u(\cdot;u_0)\in C([0,T); H^s(\R))\cap C^1([0,T); H^{s-1}(\R)).$ \
Moreover, the solution depends continuously on the initial data and is called \emph{a strong solution} due to its smoothness.
\end{prop}

Furthermore, the strong solution is a global one if the initial condition is sufficiently ``regular". More specifically, we have the following global existence result established in \cite{L-Y} and thus proof  is omitted.
\begin{prop}[\cite{L-Y}] \label{prolinfty}
Given that the initial profile $u_0\in H^s(\R)$ with $s>\frac{3}{2}$ and $w_0=u_0-u_{0,xx}+\frac{2}{3}\kappa$ is a Radon measure of fixed sign, the strong solution to the Cauchy problem \eqref{e:DP-IVP} then exists globally in time; that is,
\[
u=u(\cdot;u_0)\in C([0,\infty); H^s(\R))\cap C^1([0,\infty); H^{s-1}(\R)),
\]
which admits the following additional estimates.
\begin{enumerate}
    \item The magnitude of $u_x$ is bounded above by the sum of the magnitude of $u$ and the constant $\frac{2\kappa }{3}$. As a matter of fact, we have, for all $(t,x)\in[0,\infty)\times\R$,
    \begin{equation}\label{e:uxu}
        \vert u_x(t,x)\vert\leq\vert u(t,x)+\frac{2}{3}\kappa\vert.
    \end{equation}
    \item The $L^\infty$ norm of $u$ is bounded. More specifically, we have, for all $t\in[0,\infty)$,
    \begin{equation}\label{e:uinfty}
        \Vert u(t,\cdot)\Vert_{L^\infty}\leq  2(1+\sqrt{2})\Vert u_0\Vert_{L^2(\R)}+\frac{4}{3}\kappa.
    \end{equation}
\end{enumerate}
\end{prop}

\subsection{\emph{A priori} estimates}
To bound the high order term $\int_\R u^3dx$ in the Lagrangian with the conserved quantity $S(u)$ whose square root is equivalent to the $L^2$-norm of $u$, we take advantage of the argument in \cite[Lemma 2]{KM} to derive the following \emph{a priori} estimate.
\begin{prop}[\emph{a priori} $L^\infty$-$L^2$ estimate]\label{hlinfty}
Let ${f}\in W^{1,\infty}(\R)\cap L^2(\R)$ and the initial data $u_0\in H^s(\R)$ with $s>\frac{3}{2}$ and $w_0=m_0+\frac{2\kappa }{3}$ a Radon measure of fixed sign. The difference between the strong solution $u$ to the Cauchy problem \eqref{e:DP-IVP} and the function ${f}$, denoted as ${g}(t,x)\triangleq u(t,x)-{f}(x)$, admits the following estimate
\begin{equation}\label{e:Linfty}
\Vert {g}(t,\cdot)\Vert_{L^\infty(\R)}\leq
\Vert {g}(t,\cdot)\Vert_{L^2(\R)}^{2/3}\bigg(1+\frac{4}{3}\kappa+\sqrt{2}\Vert {g}(t,\cdot)\Vert_{L^2(\R)}^{2/3}+2\Vert {f}\Vert_{L^\infty(\R)}+2\Vert {f}'\Vert_{L^\infty(\R)}\bigg), \quad\forall t\in[0,\infty).
\end{equation}
\end{prop}
\begin{proof}
Fix $t\in[0,\infty)$, we denote $G=\Vert {g}(t,\cdot)\Vert_{L^2(\R)}^{2/3}$ and assume $G>0$, due to the fact that the case $G=0$ makes both sides of \eqref{e:Linfty} zero. Fixing $x\in\R$, there exists $k\in\mathbb{Z}$ such that
$x\in[kG,\,(k+1)G)$. By the mean value theorem, there exists
$\bar{x}\in[(k-1)G,\,kG]$ such that
\[
{g}^2(t,\bar{x})=\frac{1}{G}\int_{(k-1)G}^{kG}{g}^2(t,\eta)\rmd\eta\leq\frac{1}{G}\Vert {g}(t,\cdot)\Vert_{L^2(\R)}^2=G^2,
\]
which, together with Proposition \ref{prolinfty} and that $0\leq x-\bar{x}\leq 2G$, yields
\begin{equation}\label{ppp}
\begin{aligned}
{g}(t,x)&={g}(t,\bar{x})+\int_{\bar{x}}^x {g}_\eta(t,\eta)\rmd\eta\\
&\geq-G-\frac{4G}{3}\kappa-\sqrt{2G}\left\Vert\bigg(\vert u(t,\cdot)\vert +\vert {f}'\vert \bigg)\right\Vert_{L^2([(k-1)G,(k+1)G])}\\
&\geq-G-\frac{4G}{3}\kappa-\sqrt{2G}\left\Vert\bigg(\vert {g}(t,\cdot)\vert +\vert {f}\vert +\vert {f}'\vert \bigg)\right\Vert_{L^2([(k-1)G,(k+1)G])}\\
&\geq-G-\frac{4G}{3}\kappa-\sqrt{2G}\left[\Vert {g}(t,\cdot)\Vert_{L^2(\R)}+\sqrt{2G}\bigg(\Vert {f}\Vert_{L^\infty(\R)}+\Vert {f}'\Vert_{L^\infty(\R)}\bigg)\right]\\
&=-G\bigg(1+\frac{4}{3}\kappa+\sqrt{2}G+2\Vert {f}\Vert_{L^\infty(\R)}+2\Vert {f}'\Vert_{L^\infty(\R)}\bigg).
\end{aligned}
\end{equation}
Similarly to the argument in Proposition \ref{prolinfty}, we use the proof by contradiction and suppose that there exists $x_*\in\R$ such that
\[
{g}(t,x_*)>G\bigg(1+\frac{4}{3}\kappa+\sqrt{2}G+2\Vert {f}\Vert_{L^\infty(\R)}+2\Vert {f}'\Vert_{L^\infty(\R)}\bigg).
\]
Then there exists $k_*\in \Z$ with $x_*\in [k_*G,(k_*+1)G)$
and by the mean value theorem, there exists $\bar{x}_*\in[(k_*+1)G, (k_*+2)G]$  such that, on one hand,
\[
{g}^2(t,\bar{x}_*)=\frac{1}{G}\int_{(k_*+1)G}^{(k_*+2)G}{g}^2(t,\eta)\rmd\eta\leq G^2.
\]
On the other hand, proceeding as in \eqref{ppp}, we have
\begin{equation*}
\begin{aligned}
{g}(t,\bar{x}_*)
=&{g}(t,x_*)+\int_{x_*}^{\bar{x}_*}[u_\eta(t,\eta)-{f}'(\eta)]\rmd\eta\\
>&G\bigg(1+\frac{4}{3}\kappa+\sqrt{2}G+2\Vert {f}\Vert_{L^\infty}+2\Vert {f}'\Vert_{L^\infty}\bigg)-\frac{4G}{3}\kappa-\sqrt{2G}\left\Vert\bigg(\vert u(t,\cdot)\vert +\vert {f}'\vert \bigg)\right\Vert_{L^2([k_*G,(k_*+2)G])}\\
\geq &G\bigg(1+\frac{4}{3}\kappa+\sqrt{2}G+2\Vert {f}\Vert_{L^\infty}+2\Vert {f}'\Vert_{L^\infty}\bigg)-\frac{4G}{3}\kappa-G\left(\sqrt{2}G+2\Vert {f}\Vert_{L^\infty}+2\Vert {f}'\Vert_{L^\infty}\right)\\
=&G.
\end{aligned}
\end{equation*}
Again, the incompatibility of the above two estimates concludes the proof of the proposition.
\end{proof}

\subsection{Positivity of quadratic form}
For any fixed $c>2\kappa >0$, the second variational derivative of the Lagrangian $Q_c(u)= H(u)+cS(u)$, at the unique solitary wave $\phi=\phi_c$, admits the following expression,
\[
\cal \triangleq \frac{\delta^2 Q_c}{\delta u^2}(\phi)=-\phi_c-2\kappa (4-\partial_x^2)^{-1}+c(1-\partial_x^2)(4-\partial_x^2)^{-1}\ :  L^2(\R)\to L^2(\R),
\]
which is a well-defined, self-adjoint, bounded linear operator.
In \cite{LLW1} we proved the following spectral properties  about the operator $\cal$.
\begin{prop}\label{spectrum}
The spectrum of the operator $\cal$, denoted as $\sigma(\cal)$, admits the following properties.
\begin{enumerate}
\item $\sigma(\cal)$ lies on the real line; that is, $\sigma(\cal)\subset \R$.
\item $0\in \sigma(\cal)$ is a simple eigenvalue with $\partial_x \phi_c  $ as its eigenfunction.
\item There exists $\underline{\lambda}>0$ such that $\sigma(\cal)\cap (-\infty,0)=\{-\underline{\lambda}\}$. Moreover, $-\underline{\lambda}$ is a simple eigenvalue whose corresponding normalized eigenfunction is denoted as $\chi$.
\item The set of essential spectrum, denoted as $\sigma_{{\rm ess}}(\cal)$, lies on the positive real axis, admitting a positive distance to the origin.
\end{enumerate}
\end{prop}

{Based on the spectral decomposition of $\cal$, we have
\begin{equation}\label{e:di-sum}
L^2(\R)=\spn\{\chi\}\oplus \spn\{\partial_x\phi\}\oplus X_+,\quad
\end{equation}
where $X_+\triangleq(\spn\{\chi, \partial_x\phi\})^\perp$ is invariant under $\cal$ and there exists $\overline{\lambda}>0$ such that
\[
(\cal p,p)\geq\overline{\lambda}\Vert p\Vert_{L^2(\R)},\quad \mbox{ for any } p\in X_+.
\]
Moreover, we have the following lemma.}
\begin{lemma} \label{Positivity}
There exists $\eta>0$ sufficiently small, such that for any $y\in L^2(\R)$ satisfying
\begin{equation}\label{orth}
\left\vert(y,\phi)_{\cas}\right\vert+\left\vert(y,\partial_x\phi)_{\cas}\right\vert \leq \eta\Vert y\Vert_{L^2(\R)},
\end{equation}
there exists $\theta>0$ such that $(\cal y,y)\geq\theta\Vert y\Vert_{L^2(\R)}^2.$
\end{lemma}

\begin{proof}
 The proof consists of two steps: first we deal with the case  $\eta=0$  then the case  $\eta>0$. In order to do that, we first develop a more friendly expression of $(y,\phi)_{\cas}$.
Recalling  the soliton equation
$\frac{\delta H}{\delta u}(\phi)+c\frac{\delta S}{\delta u}(\phi)=0,$
which, taken derivative with respect to $c$, yields $
\cal\partial_c\phi=-\frac{\delta S}{\delta u}(\phi).$
We have
\[
(y,\phi)_{\cas}=(y,\frac{\delta S}{\delta u}(\phi))=-(y,\cal\partial_c\phi).
\]
Furthermore, we have
\[
(\cal\partial_c\phi,\partial_c\phi)=(-\frac{\delta S}{\delta u}(\phi),\partial_c\phi)
=-\frac{\rmd}{\rmd c}S(\phi)<0,
\]
where we refer to Lemma \ref{quantity} for the last inequality. It follows from the direct sum \eqref{e:di-sum} that there exist unique $a_0, b_0\in\R$, $p_0\in X_+$, $a_0\neq0$, such that
\begin{equation}\label{pc}
\partial_c\phi=a_0(\chi+p_0)+b_0\partial_x\phi.
\end{equation}
We then have
\[
(\cal\partial_c\phi,\partial_c\phi)=a_0^2\big[(\cal \chi,\chi)+(\cal p_0,p_0)\big]=a_0^2\big[(\cal p_0,p_0)-\underline{\lambda}\big]<0,
\]
which indicates that there exists $\delta>0$ such that
\[(\cal p_0,p_0)<\textcolor{red}{-}\underline{\lambda}.\]
We can take a small positive number $\delta>0$ such that
\[
(\cal p_0,p_0)<\frac{1}{1+\delta}\underline{\lambda}.
\]
As a result, we summarize that
\begin{equation}\label{e:yphi}
(y,\phi)_{\cas}=a_0\big[ \underline{\lambda}(y,\chi)-(y,\cal p_0)\big].
\end{equation}
We now ready to prove the theorem for the case $\eta=0$.\\ \vskip 0.01in
\noindent {\it Case 1 ($\eta=0$).} If $y_1\in L^2(\R)$ satisfies \eqref{orth} for $\eta=0$, then there exists $\theta_1>0$ such that
\begin{equation}\label{step1}
(\cal y_1,y_1)\geq\theta_1\Vert y_1\Vert_{L^2(\R)}^2.
\end{equation}
\noindent {\it Case 1.1.} If
\[y_1=b\partial_x\phi+\tilde{p}_1,\quad b\in\R, \quad\tilde{p}_1\in X_+,\]
then
\[
\begin{aligned}
\Vert\tilde{p}_1\Vert_{L^2(\R)}^2
=&\Vert y_1-b\partial_x\phi\Vert_{L^2(\R)}^2
\geq(y_1-b\partial_x\phi,y_1-b\partial_x\phi)_{\cas}
=\Vert y_1\Vert_{\cas}^2+b^2\Vert\partial_x\phi\Vert_{\cas}^2\\
\geq&\Vert y_1\Vert_{\cas}^2\geq\frac{1}{4}\Vert y_1\Vert_{L^2(\R)}^2.
\end{aligned}
\]
Therefore, we have
\begin{equation}\label{ss1}
(\cal y_1,y_1)=(\cal \tilde{p}_1,\tilde{p}_1)\geq \overline{\lambda}\Vert \tilde{p}_1\Vert_{L^2(\R)}^2\geq\frac{1}{4}\overline{\lambda}\Vert y_1\Vert_{L^2(\R)}^2.
\end{equation}
\noindent {\it Case 1.2.} If
\[y_1=b\partial_x\phi+a(\chi+p_1),\quad a,b\in\R, \quad a\neq 0,\quad p_1\in X_+,\]
then as in {\it Case 1.1}, we get
\begin{equation}\label{chip}
a^2\Vert\chi+p_1\Vert_{L^2(\R)}^2
\geq\frac{1}{4}\Vert y_1\Vert_{L^2(\R)}^2.
\end{equation}
In addition, by \eqref{orth} and \eqref{e:yphi}, we have
\[
\begin{aligned}
0
&=\vert(y_1,\phi)_{\cas}\vert=\vert a_0\big[ \underline{\lambda}(y_1,\chi)-(y_1,\cal p_0)\big] \vert
&=\vert a_0a\vert\cdot\vert(\chi+p_1,\cal(\chi+p_0))\vert
=\vert a_0a\vert \cdot\vert \underline{\lambda}-(\cal p_0,p_1)\vert,
\end{aligned}
\]
yielding
\begin{eqnarray*}
\cal p_0,p_1)=\underline{\lambda}.
\end{eqnarray*}
As a result, we derive two different lower bounds of $(\cal(\chi+p_1),\chi+p_1)$; that is,
\[
\begin{aligned}
(\cal(\chi+p_1),\chi+p_1)=&(\cal p_1,p_1)-\underline{\lambda}\geq\frac{(\cal p_0,p_1)^2}{(\cal p_0,p_0)}-\underline{\lambda}\geq\frac{\underline{\lambda}^2}{(\cal p_0,p_0)}-\underline{\lambda}
\geq\frac{\underline{\lambda}^2}{\frac{1}{1+\delta}\underline{\lambda}}-\underline{\lambda}=\delta\underline{\lambda},\\
(\cal(\chi+p_1),\chi+p_1)=&(\cal p_1,p_1)-\underline{\lambda}\geq\overline{\lambda}\Vert p_1\Vert_{L^2(\R)}^2-\underline{\lambda}.
\end{aligned}
\]
which, in turns, yield
\begin{equation}\label{ss2}
\begin{aligned}
(\cal y_1,y_1)=&a^2(\cal(\chi+p_1),\chi+p_1)=a^2\frac{(\delta+2)(\cal(\chi+p_1),\chi+p_1)}{\delta+2}\geq a^2\frac{\delta[\overline{\lambda}\Vert p_1\Vert_{L^2(\R)}^2-\underline{\lambda}]+2\delta\underline{\lambda}}{\delta+2}
\\
\geq&a^2\frac{\delta\min\{\underline{\lambda},\overline{\lambda}\}}{\delta+2}(\Vert p_1\Vert_{L^2(\R)}^2+1)=a^2\frac{\delta\min\{\underline{\lambda},\overline{\lambda}\}}{\delta+2}\Vert p_1+\chi\Vert_{L^2(\R)}^2\geq \frac{\delta\min\{\underline{\lambda},\overline{\lambda}\}}{4(\delta+2)}\Vert y_1\Vert_{L^2(\R)}^2.
\end{aligned}
\end{equation}
Case 1 then follows from \eqref{ss1} and \eqref{ss2}.

\noindent {\it Case 2 ($\eta>0$). }  If $y\in L^2(\R)$ satisfies \eqref{orth} for $0<\eta\ll 1$, then there exists $\theta>0$ such that
\begin{equation}\label{step1}
(\cal y,y)\geq\frac{1}{4}\theta\Vert y\Vert_{L^2(\R)}^2.
\end{equation}
The function $y$ admits the decomposition
\begin{equation}\label{orthoo}
y=y_1+a_1\phi+b_1\partial_x\phi, \quad\mbox{with} \quad (y_1,\phi)_{\cas}=(y_1,\partial_x\phi)_{\cas}=0,
\end{equation}
where $a_1=\frac{(y,\phi)_{\cas}}{(\phi,\phi)_{\cas}},\quad b_1=\frac{(y,\partial_x\phi)_{\cas}}{(\partial_x\phi,\partial_x\phi)_{\cas}},\quad\vert a_1\vert\leq\frac{4\eta}{\Vert\phi\Vert_\cas^2}\Vert y\Vert_{\cas},\quad\vert b_1\vert\leq\frac{4\eta}{\Vert\partial_x\phi\Vert_\cas^2}\Vert y\Vert_{\cas},$
yielding
\[
\begin{aligned}
\Vert y_1\Vert_{\cas}^2=\Vert y\Vert_{\cas}^2- \Vert a_1\phi+b_1\partial_x\phi \Vert_{\cas}^2\geq &\Vert y\Vert_{\cas}^2-(|a_1|\Vert \phi\Vert_{\cas}+|b_1|\Vert \partial_x\phi\Vert_{\cas})^2\\
\geq & \left[ 1-\left( \frac{4\eta}{\min\{ \Vert \phi\Vert_{\cas},  \Vert \partial_x\phi\Vert_{\cas}\}}\right)^2\right]\Vert y\Vert_{\cas}^2.
\end{aligned}
\]
Taking $\eta\leq \frac{\sqrt{3}}{8}\min\{ \Vert \phi\Vert_{\cas},  \Vert \partial_x\phi\Vert_{\cas}\},$
we have
\begin{equation}
\Vert y_1\Vert_{L^2(\R)}\geq\Vert y_1\Vert_{\cas}\geq\frac{1}{2}\Vert y\Vert_{\cas}\geq \frac{1}{8}\Vert y\Vert_{L^2(\R)}.
\end{equation}
As a result, we conclude that
\[
\begin{aligned}
(\cal y,y)
=&(\cal y_1,y_1)+(\cal(a_1\phi+b_1\partial\phi),a_1\phi+b_1\partial\phi)+2(\cal(a_1\phi+b_1\partial\phi),y_1)\\
=&(\cal y_1,y_1)+a_1^2(\cal \phi,\phi)+2a_1(\cal\phi,y_1)\\
\geq &\frac{\delta\min\{\underline{\lambda},\overline{\lambda}\}}{4(\delta+2)}\Vert y_1\Vert_{L^2(\R)}^2-\frac{\eta^2\vert (\cal\phi,\phi)\vert}{\Vert \phi\Vert_{\cas}^4}\Vert y\Vert_{L^2(\R)}^2-\frac{2\eta\Vert \cal \phi\Vert_{L^2(\R)}}{\Vert\phi\Vert_{\cas}^2}\Vert y\Vert_{L^2(\R)}\Vert y_1\Vert_{L^2(\R)}\\
\geq &\underbrace{\left[ \frac{\delta\min\{\underline{\lambda},\overline{\lambda}\}}{32(\delta+2)}-\frac{\Vert\cal\phi\Vert_{L^2(\R)}\Vert\phi\Vert_{L^2(\R)}}{\Vert \phi\Vert_{\cas}^4}\eta^2-\frac{8\Vert \cal \phi\Vert_{L^2(\R)}}{\Vert\phi\Vert_{\cas}^2}\eta\right] }_{:=f(\eta)}\Vert y\Vert_{L^2(\R)}^2
\end{aligned}
\]
It is straightforward to see that $f(\eta)=0$ admits two root $\eta_-<0<\eta_+$. Taking
\[
0<\eta<\min\{\frac{\eta_+}{2}, \frac{\sqrt{3}}{8}\min\{ \Vert \phi\Vert_{\cas},  \Vert \partial_x\phi\Vert_{\cas}\}\}, \quad \theta=f(\frac{\eta_+}{2})>0,
\]
we have $(\cal y,y)\geq \theta\Vert y\Vert_{L^2(\R)}^2,$
which concludes the proof.
\end{proof}

\begin{lemma}\label{quantity}
For the $c$-speed solitary wave $\phi=\phi(x; c)$, we have
\begin{equation}\label{derivative}
\frac{\rmd}{\rmd c}H(\phi )=-c\frac{\rmd}{\rmd c}S(\phi )=-\frac{3c^2(c+\kappa)}{2(3c+2\kappa )^2}\sqrt{c^2-2c\kappa }<0.
\end{equation}
\end{lemma}

\begin{proof}
The second equality was proved in \cite{LLW1} and we only prove the first equality. Note that $\phi(x;c)$ is the critical point of the Lagrangian $Q_{\lambda}(u)=H(u)+\lambda S(u)$ with $\lambda=c$; that is,
\begin{equation}\label{e:soliton-eq}
\left(\frac{\delta H}{\delta u}(\phi)+\lambda\frac{\delta S}{\delta u}(\phi)\right)\bigg\vert_{\lambda=c}=0,
\end{equation}
we have
\[
\frac{\rmd}{\rmd c}H(\phi )+c\frac{\rmd}{\rmd c}S(\phi )=\frac{\rmd Q_{\lambda}}{\rmd c}(\phi)\bigg\vert_{\lambda=c}= \left < \left(\frac{\delta H}{\delta u}(\phi)+\lambda\frac{\delta S}{\delta u}(\phi)\right)\bigg\vert_{\lambda=c}, \, \frac{\partial \phi}{\partial c} \right > = 0.
\]
\end{proof}
\section{Decomposition and properties of solution near N solitons}

Fix a constant $\sigma_0$ such that $0<\sigma_0<\frac{1}{2}\min\left\{\sqrt{1-\frac{2\ka}{c_1^0}}, c_1^0-2\ka, c_2^0-c_1^0, c_3^0-c_2^0,\cdots, c_N^0-c_{N-1}^0\right\}$. Let $\gamma_0=\min\{1/(8B), \sigma_0/8\}$, where $B$ is the constant determined in Lemma \ref{Posi}. We have by \eqref{exp0} that
\begin{equation}\label{exp1}
\vert\phi_c(\cdot)\vert\leq Ce^{-2\sigma_0\vert\cdot\vert}.
\end{equation}
On the other hand, there appears 
\begin{equation}\label{exp2}
\vert\left((1-\partial_x^2)^{-1}\phi_c(\cdot-m)\right)(\xi)\vert=\frac{1}{2}\aaa e^{-\vert \xi-r\vert}\phi_c(r-m) dr\leq C\aaa e^{-\vert \xi-r\vert-2\sigma_0\vert r-m\vert}dr\leq Ce^{-\sigma_0\vert \xi-m\vert},
\end{equation}
\begin{equation}\label{exp3}
\vert\left((4-\partial_x^2)^{-1}\phi_c(\cdot-m)\right)(\xi)\vert=\frac{1}{4}\aaa e^{-2\vert \xi-r\vert}\phi_c(r-m) dr \leq Ce^{-\sigma_0\vert \xi-m\vert}.
\end{equation}
For  $A_0,\ L,\ \alpha>0$, we define
\begin{align*}
\cau(L,\al)\triangleq&\left\{
u\in L^2(\R)\ : \inf_{x_j-x_{j-1}\geq L}\left\Vert u-\sum_{j=1}^N \phi_{c_j^0}(\cdot-x_j)\right\Vert_{L^2}\leq \al
\right\},\\
\cav_{A_0}(L,\al)\triangleq&\left\{
u\in L^2(\R)\ : \inf_{x_j-x_{j-1}\geq L}\left\Vert u-\sum_{j=1}^N \phi_{c_j^0}(\cdot-x_j)\right\Vert_{L^2}\leq A_0(\al+e^{-\gamma_0L/2})
\right\}.
\end{align*}

Our strategy is to prove that there exists $A_0>0$, $L_0>0$, and $\al_0>0$ such that, $\forall u_0\in H^s(\R)$, if for some $L>L_0$, $\al<\al_0$, $\Vert u_0-\sum_{j=1}^N\phi_{c_j^0}(\cdot-x_j^0)\Vert_{L^2}\leq\al$, where $x_j^0>x_{j-1}^0+L$, then $\forall t\geq 0$, $u(t)\in\cav_{A_0}(L,\al)$. In order to establish Theorem \ref{Nstability}, it suffices to prove the following
\begin{prop} \label{priori}
There exists $A_0>0$, $L_0>0$, and $\al_0>0$ such that, for all $u_0\in H^s(\R)$, if
\begin{equation}\label{initial}
\left\Vert u_0-\sum_{j=1}^N\phi_{c_j^0}(\cdot-x_j^0)\right\Vert_{L^2}\leq\al,
\end{equation}
where $0<\al<\al_0$, $x_j^0>x_{j-1}^0+L$, $L>L_0$, and if for $t^*>0$,
\begin{equation}\label{assum}
\forall t\in[0,t^*], \quad u(t)\in\cav_{A_0}(L,\al),
\end{equation}
then
\begin{equation}
\forall t\in[0,t^*], \quad u(t)\in\cav_{A_0/2}(L,\al),
\end{equation}
where $A_0,\ L_0, $ and $ \al$ are independent of $t^*$.
\end{prop}
The proof of Proposition \ref{priori} is approached via a series of lemmas, the proofs of which are provided in the next section.
\begin{lemma}[{Decomposition of the solution}]\label{Decom}
There exist $L_1,\ \al_1, \ K_1>0$ such that the following is true: If for some $t^*>0$, the solution of \eqref{DP} $u(t,x)\in \cau(L,\al)$ for all $t\in[0,t^*]$ with
$L>L_1,\ 0<\al<\al_1$,
then there exist unique $C^1$ functions $c_j\ :\  [0,t^*]\to(2\ka,\infty),\ x_j\ :\ [0,t^*]\to\R$, such that
$$\e(t,x)\triangleq u(t,x)-\sum_{j=1}^NR_j(t,x),\quad\mbox{ where }\quad R_j(t,x)\triangleq \phi_{c_j(t)}(x-x_j(t)),$$
satisfies the following orthogonality conditions:
\begin{equation}\label{ortho0}
\forall j, \forall t\in[0,t^*],\quad \int (1-\partial_x^2)R_j(t)\widehat{\e}(t)dx=\int (1-\partial_x^2)R_{j,x}(t)\widehat{\e}(t)dx=0.
\end{equation}
Moreover, for any $t\in[0,t^*]$,
\begin{equation}\label{c1-0}
\Vert \e(t,\cdot)\Vert_{L^2}+\sum_{j=1}^N\vert c_j(t)-c_j^0\vert \leq K_1\al, \quad x_j(t)>x_{j-1}(t)+L-K_1\al.
\end{equation}
\begin{equation}\label{c2-0}
\vert \dot{c}_j(t)\vert+\vert \dot{x}_j(t)-c_j(t)\vert\leq K_1\left(\int  \e^2 e^{- \sigma_0\vert x-x_j(t)\vert} dx\right)^{1/2}+K_1e^{-\sigma_0L/2}.
\end{equation}
\end{lemma}

\begin{lemma}[{Energy bound}]\label{Energy}
Then there exist $K_2>0$ and $L_2>0$ such that if $L_0>L_2$,  then, for all $t\in[0,t^*]$:
\begin{equation*}
\left\vert
\sum_{j=1}^N\left[H(\ppt)-H(\ppo)\right]-\int \left[\frac{1}{2}R(t)\e^2(t)+\kappa\e(t)\widehat{\e}(t)\right] dx
\right\vert
\end{equation*}
\begin{equation}\label{16}
\leq K_2\left(\Vert\e(t)\Vert_{L^3}^3+\Vert\e(0)\Vert_{L^2}^2+\Vert\e(0)\Vert_{L^3}^3+e^{-\frac{\sigma_0L}{2}}\right).
\end{equation}
\end{lemma}

Let $\psi(x)=\frac{2}{\pi}\arctan(e^\frac{x}{B})$ and define for $j\geq 2$, $m_j(t)\triangleq\frac{x_{j-1}(t)+x_j(t)}{2}$ and
\begin{equation}
I_j(t)\triangleq\frac{1}{2}\int \psi(x-m_j(t))\left[4\widehat{u}^2(t)+5\widehat{u}_x^2(t)+\widehat{u}_{xx}^2(t)\right] dx,\quad\quad
\end{equation}
where $B$ is a large constant chosen in the proof of Lemma \ref{Posi}.
Note that $\psi\in(0,1)$ and $\psi$ is monotonically increasing, satisfying
\begin{equation}\label{psi}
\vert\psi^{(4)}\vert\leq \frac{3}{B^3}\psi',\quad\vert\psi'''\vert\leq\frac{1}{B^2}\psi',\quad\psi''\leq\frac{1}{B}\psi',\quad\mbox{and }\ \psi^{(k)}(x)\leq Ce^{-\frac{\vert x\vert}{B}}\mbox{ for } k=1,2,3,4.
\end{equation}

\begin{lemma}[{Monotonicity}]\label{monotone}
There exist $K_3>0,\ L_3>0,\  \al_2>0$ such that if $L_0>L_3$, $\al_0<\al_2$,
then for all $t\in[0,t^*]$:
\begin{equation*}
I_j(t)-I_j(0)\leq K_3e^{-\frac{L}{4B}}.
\end{equation*}
\end{lemma}

Let $\al_0$ and $L_0$ be chosen such that \eqref{co0}, \eqref{co1}, and \eqref{co2} be satisfied. Then we have

\begin{lemma}[{Positivity}] \label{Posi} There exists $L_4>0$ and $\lambda_0>0$ such that if 
$L_0>L_4$, then, for all $t\in[0,t^*]$,
\begin{equation}\label{positivity}
\int \left[-R(t)\e^2-2\kappa \e\widehat{\e}+c(t,x)(4\widehat{\e}^2+5\widehat{\e}_x^2+\widehat{\e}_{xx}^2)\right]dx\geq \lambda_0\Vert\e(t)\Vert_{L^2}^2,
\end{equation}
where $c(t,x)\triangleq c_1(t)+\sum_{j=2}^{N}[c_j(t)-c_{j-1}(t)]\psi(x-m_j(t))$.
\end{lemma}

\begin{lemma}[{ Quadratic control}]\label{Qua} There exists $K_4>0$ such that, for any $t\in[0,t^*]$,
\begin{equation}\label{quadratic}
\sum_{j=1}^N\vert c_j(t)-c_j(0)\vert\leq K_4(\Vert\e(t)\Vert_{L^2}^2+\Vert\e(0)\Vert_{L^2}^2+e^{-\frac{1}{B}(\frac{\sigma_0t}{2}+\frac{L}{8})}+e^{-\sigma_0L/8}).
\end{equation}
\end{lemma}

\begin{lemma} \label{Lya}There exists $K_5>0$ independent of $A_0$, such that, for any $t\in[0,t^*]$,
\begin{equation}\label{lya}
\Vert\e(t)\Vert_{L^2}^2\leq K_5\left(\Vert\e(0)\Vert_{L^2}^2+e^{-\gamma_0L}\right).
\end{equation}
\end{lemma}
With these lemmas in hand, we are in a position to prove Proposition \ref{priori}.
\begin{proof}[Proof of Proposition \ref{priori}]  Let $A_0>1$  be fixed later. We choose $\al_0$ and $L_0(\geq L_1)$ such that
\begin{equation}\label{co0}
A_0\left(\al_0+ e^{-\gamma_0L_0/2}\right)<\al_1,
\end{equation}
where $\al_1$ and $L_1$ are defined in Lemma \ref{Decom}. Therefore, assuming \eqref{assum}, Lemma \ref{Decom} is applicable: there exist unique $C^1$ functions $c_j :  [0,t^*]\to(2\ka,\infty),\ x_j : [0,t^*]\to\R$, such that
$$\e(t,x)\triangleq u(t,x)-\sum_{j=1}^NR_j(t,x)\triangleq u(t,x)-R(t,x),\quad\mbox{ where }\quad R_j(t,x)=\phi_{c_j(t)}(x-x_j(t)),$$
satisfies the following orthogonality conditions:
\begin{equation}\label{ortho}
\forall j, \forall t\in[0,t^*],\quad \int (1-\partial_x^2)R_j(t)(4-\partial_x^2)^{-1}\e(t)dx=\int (1-\partial_x^2)R_{j,x}(t)(4-\partial_x^2)^{-1}\e(t)dx=0.
\end{equation}
Moreover, for any $t\in[0,t^*]$,
\begin{equation}\label{c1}
\Vert \e(t,\cdot)\Vert_{L^2}+\sum_{j=1}^N\vert c_j(t)-c_j^0\vert \leq K_1A_0(\al+e^{-\gamma_0L/2}), \quad x_j(t)>x_{j-1}(t)+L-K_1A_0(\al+e^{-\gamma_0L/2}),
\end{equation}
\begin{equation}\label{c2}
\vert \dot{c}_j(t)\vert+\vert \dot{x}_j(t)-c_j(t)\vert\leq K_1\left(\int \e^2 e^{-\sigma_0\vert x-x_j(t)\vert} dx\right)^{1/2}+K_1e^{-\gamma_0L/2}.
\end{equation}
Note also that by \eqref{initial} and \eqref{000}, we have stronger estimates for initials as follows
\begin{equation}\label{c3}
\Vert \e(0)\Vert_{L^2}+\sum_{j=1}^N\vert c_j(0)-c_j^0\vert \leq K_1 \al , \quad x_j(0)>x_{j-1}(0)+L-K_1\al.
\end{equation}
From \eqref{c1}, \eqref{c2}, and \eqref{c3}, we can choose $\al_0$ small enough, and $L_0$ large enough, such that $\forall t\in[0,t^*]$,
\begin{equation}\label{co1}
\dot{x}_1(t)\geq2\ka+\sigma_0,\quad \dot{x}_j(t)-\dot{x}_{j-1}(t)\geq \sigma_0,\quad c_1(t)>2\ka+\sigma_0, \quad c_j(t)-c_{j-1}(t)\geq\sigma_0,
\end{equation}
\begin{equation}\label{co2}
x_j(t)-x_{j-1}(t)\geq L/2>L_0/2,\quad \Vert\e(t)\Vert_{L^2}\leq \sigma_0/8.
\end{equation}

In view of  \eqref{c1}, \eqref{c3}, \eqref{quadratic}, and \eqref{lya}, we have
\begin{eqnarray*}
&&\left\Vert u(t)-\sum \phi_{c_j^0}(\cdot-x_j(t))\right\Vert_{L^2}\\
&\leq&\left\Vert \e(t)\right\Vert_{L^2}+C\sum \vert c_j(t)-c_j(0)\vert+C\sum \vert c_j(0)-c_j^0\vert\\
&\leq&\left\Vert \e(t)\right\Vert_{L^2}+CK_4(\Vert\e(t)\Vert_{L^2}^2+\Vert\e(0)\Vert_{L^2}^2+e^{-\gamma_0L})+CK_1\al\\
&\leq&\sqrt{K_5(\Vert\e(0)\Vert_{L^2}^2+e^{-\gamma_0L})}+CK_4(1+K_5)(\Vert\e(0)\Vert_{L^2}^2+e^{-\gamma_0L})+CK_1\al\\
&\leq&K_6(\al+e^{-\gamma_0L/2}),
\end{eqnarray*}
for some constant $K_6>0$.
Choosing $A_0=4K_6$, this completes the proof of Proposition \ref{priori}.
\end{proof}
\section{Proof of Lemmas}
We shall make use of the following essentially
\begin{equation}\label{eeee}
\Vert \e(t,\cdot)\Vert_{L^\infty}\leq C\Vert \e(t,\cdot)\Vert_{L^2}^{\frac{2}{3}}\leq C,
\end{equation}
which follows from applying \eqref{e:Linfty} for $g=\e(t,x)=u-R$.
\begin{proof}[Proof of Lemma \ref{Decom}] We make the following claim:

{\it \begin{itemize}\item
There exists $L_1,\ \al_1>0$ and unique $C^1$ functions $(c_j,x_j)\ : \cau(L_1,\al_1)\to (2\ka,\infty)\times\R$, where $j=1, 2,\cdots, N$, such that if $u\in\cau(L_1,\al_1)$, and
\[
\e(x)\triangleq u(x)-\sum_{j=1}^N\phi_{c_j(u)}(x-x_j(u)),
\]
then for any $j=1,\cdots, N$,
\begin{equation*}
\quad \int (1-\partial_x^2)\phi_{c_j(u)}(x-x_j(u))\widehat{\e}(x)dx=\int (1-\partial_x^2)\phi_{c_j(u),x}(x-x_j(u))\widehat{\e}(x)dx=0.
\end{equation*}
Moreover, There exists $K_1>0$ such that if $u\in\cau(L,\al)$, with $0<\al<\al_1,\ L>L_1$, then
\begin{equation}\label{000}
\Vert \e\Vert_{L^2}+\sum_{j=1}^N\vert c_j(u)-c_j^0\vert \leq K_1\al,\quad x_j(u)>x_{j-1}(u)+L-K_1\al.
\end{equation}\end{itemize}
}
This is a standard application of the implicit function theorem.  The non-degeneracy conditions in applying the implicit function theorem are
\begin{eqnarray*}
&&\int (1-\partial_x^2)\phi_{c_j}(x-r_j^0)\frac{d}{dc_j}(4-\partial_x^2)^{-1}\phi_{c_j}(x-r_j^0)\ dx=\frac{d}{dc_j}S(\phi_{c_j})\overset{\eqref{derivative}}{>}0,\\
&&\int (1-\partial_x^2)\partial_x\phi_{c_j}(x-r_j^0)\frac{d}{dc_j}(4-\partial_x^2)^{-1}\phi_{c_j}(x-r_j^0)\ dx\triangleq b_j,\\
&&\int (1-\partial_x^2)\phi_{c_j}(x-r_j^0)\frac{d}{dx}(4-\partial_x^2)^{-1}\phi_{c_j}(x-r_j^0)\ dx =0,\\
&&\int (1-\partial_x^2)\partial_x\phi_{c_j}(x-r_j^0)\frac{d}{dx}(4-\partial_x^2)^{-1}\phi_{c_j}(x-r_j^0)\ dx= 2S(\partial_x\phi_{c_j^0})>0;
\end{eqnarray*}
and for $k\neq j$,
\begin{eqnarray}\label{mix}
\nonumber&&\left\vert-\int (1-\partial_x^2)\phi_{c_j}(x-r_j^0)\frac{d}{dc_k}(4-\partial_x^2)^{-1}\phi_{c_k}(x-r_k^0)\ dx\right\vert
\\
&\leq& C\int e^{-\sigma_0(2\vert x-r_j^0\vert+\vert x-r_k^0\vert)}dx\leq Ce^{-\sigma_0\vert r_j^0-r_k^0\vert}\overset{\eqref{exp3}}{\leq} Ce^{-\sigma_0L/2}.
\end{eqnarray}
We omit the details. (See for example \cite{MMT}, Lemma 1 for instance.)

To continue the proof of Lemma \ref{Decom}, let $c_j(t)\triangleq c_j(u(t,\cdot)),\ x_j(t)\triangleq x_j(u(t,\cdot))$, where with some abuse of notation, the left-sided $c_j,\ x_j$ are functions of time $t$ and the right-sided $c_j,\ x_j$ are functions of $u\in\cau(L,\al)$ and are constructed as above.  Then the orthogonality conditions and \eqref{c1-0} are clearly satisfied.

We prove \eqref{c2-0} next. Applying $(4-\partial_x^2)$ on both sides of the soliton profile equation
\begin{equation}\label{profile}
c(1-\partial_\xi^2)(4-\partial_\xi^2)^{-1}\phi-\left[\frac{1}{2}(\phi)^2+(4-\partial_\xi^2)^{-1}2\kappa \phi\right]=0
\end{equation}
to get $2\ka R_{j,x}=\left[c_j(t)(1-\partial_x^2)R_j+R_jR_{j,xx}+ R_{j,x}^2-2R_j^2\right]_x.$
Therefore, from the equation \eqref{CH}:
\begin{align*}
\left(uu_{xx}+u_x^2-2u^2\right)_x
=&
(1-\partial_x^2)u_t+2\ka u_x \\
=&
(1-\partial_x^2)\left(\e_t+\frac{d}{dt}\sum R_j\right)+2\ka \left(\e_x+\sum R_{j,x}\right)\\
=&
(1-\partial_x^2)\e_t +(1-\partial_x^2)\sum \left(\dot{c}_j(t)\frac{d}{dc}R_j-\dot{x}_j(t)R_{j,x}\right)+2\ka  \e_x+2\ka\sum R_{j,x}\\
=&
(1-\partial_x^2)\e_t +(1-\partial_x^2)\sum \left(\dot{c}_j(t)\frac{d}{dc}R_j-\dot{x}_j(t)R_{j,x}\right)+2\ka  \e_x \\
&+\sum\left[c_j(t)(1-\partial_x^2)R_j+R_jR_{j,xx}+ R_{j,x}^2-2R_j^2\right]_x,
\end{align*}
which yields
\begin{align*}
(1-\partial_x^2)\e_t +2\ka  \e_x&+(1-\partial_x^2)\sum \left(\dot{c}_j(t)\frac{d}{dc}R_j -\dot{x}_j(t)R_{j,x}\right)+\sum c_j(t)(1-\partial_x^2)R_{j,x}\\
=&
\left(uu_{xx}+u_x^2-2u^2 -\sum\left(R_j R_{j,xx}+R_{j,x}^2-2R_j^2\right) \right)_x\\
=&\left[-(4-\partial_x^2)(\frac{1}{2}u^2)+(4-\partial_x^2)(\frac{1}{2}\sum R_j^2)\right]_x\\
=&\left[-(4-\partial_x^2)(\frac{1}{2}(\e+\sum R_j)^2)+(4-\partial_x^2)(\frac{1}{2}\sum R_j^2)\right]_x\\
=&-(4-\partial_x^2)(\frac{1}{2}\e^2)_x-\sum(4-\partial_x^2)(\e R_j)_x-\sum_{j\neq k}(4-\partial_x^2)(R_jR_k)_x.
\end{align*}
Taking the $L^2$ scalar product of the above equation with $(4-\partial_x^2)^{-1}R_j$, integrating by parts, using the orthogonality conditions \eqref{ortho} as well as  arguments of \eqref{mix} for mixed terms, and noting that
\begin{align*}
\int (1-\partial_x^2)\e_t(4-\partial_x)^{-1}R_jdx
=&\frac{d}{dt}\int (1-\partial_x^2)\e(4-\partial_x)^{-1}R_jdx-\int (1-\partial_x^2)\e\frac{d}{dt}(4-\partial_x)^{-1}R_jdx\\
=&-\dot{c}_j(t)\int (1-\partial_x^2)\e(4-\partial_x)^{-1}\frac{d}{dc}R_j dx=O\left(\left\Vert\e(t)\sqrt{\left\vert\frac{d}{dc}R_j\right\vert}\right\Vert_{L^2}\right),
\end{align*}
which is small for $\al$ small, we get
$$\left[\frac{d}{dc}S(\ppt)+O\left(\left\Vert\e(t)\sqrt{\left\vert\frac{d}{dc}R_j\right\vert}\right\Vert_{L^2}\right)\right]\vert\dot{c}_j(t)\vert=
O\left(\int \e^2e^{-\sigma_0\vert x-x_j(t)\vert} dx\right)^{1/2}+O(e^{-\sigma_0L/2}).
$$
Therefore, $\vert \dot{c}_j(t)\vert \leq K_1\left(\int \e^2e^{-\sigma_0\vert x-x_j(t)\vert} dx\right)^{1/2}+K_1e^{-\sigma_0L/2}.$
Taking now the scalar product as above with $(4-\partial_x^2)^{-1}R_{j,x}$ instead, arguing as above and using as well the last inequality for $\vert \dot{c}_j(t)\vert$, we get $\vert c_j(t)-\dot{x}_j(t)\vert \leq K_1\left(\int  \e^2 e^{-\sigma_0\vert x-x_j(t)\vert} dx\right)^{1/2}+K_1e^{-\sigma_0L/2}.$  So \eqref{c2} follows.
\end{proof}

\begin{proof}[Proof of Lemma \ref{Energy}] We have by \eqref{profile} and \eqref{ortho0}:
\begin{eqnarray*}
-H(u(t))
&=&\frac{1}{6}\int \left[(\e+R)^3+6\kappa (\e+R)(\widehat{\e}+\widehat{R})\right]dx\\
&=&-H(R(t))+\int \left(\frac{1}{2}R^2\e+2\kappa \widehat{R}\e\right)dx+\int \left(\frac{1}{2}\e^2R+\kappa\e\widehat{\e}\right)dx+\frac{1}{6}\int \e^3dx\\
&=&-\sum_{j=1}^NH(R_j(t))+\int \sum_{j=1}^N\left(\frac{1}{2}R_j^2\e+2\kappa \widehat{R_j}\e\right)dx+O\left(e^{-\frac{\sigma_0L}{2}}\right)\\
&&\quad\quad\quad\quad\quad+\int \left(\frac{1}{2}\e^2R+\kappa\e\widehat{\e}\right)dx+\frac{1}{6}\int \e^3dx\quad(\mbox{ by arguing as for } \eqref{mix})\\
&=&-\sum_{j=1}^NH(\phi_{c_j(t)})+\int \sum_{j=1}^N\left(c_j(t)(1-\partial_x^2)\widehat{R_j}\right)\e dx+O\left(e^{-\frac{\sigma_0L}{2}}\right)\\
&&\quad\quad\quad\quad\quad+\int \left(\frac{1}{2}\e^2R+\kappa\e\widehat{\e}\right)dx+\frac{1}{6}\int \e^3dx \\
&=&-\sum_{j=1}^NH(\phi_{c_j(t)})+O\left(e^{-\frac{\sigma_0L}{2}}\right)+\int \left(\frac{1}{2}\e^2R+\kappa\e\widehat{\e}\right)dx+\frac{1}{6}\int \e^3dx.\\
\end{eqnarray*}
Since $H(u(t))=H(u(0))$, the above implies \eqref{16}.
\end{proof}

\begin{proof}[Proof of Lemma \ref{monotone}]
Integrating by parts many times yields the following equality.
\begin{align}\label{virialpsi}
\begin{split}
\frac{d}{dt}\int (4\widehat{u}^2+&5\widehat{u}_x^2+\widehat{u}_{xx}^2)\psi(x-m_j(t)) \\
=\ &\frac{2}{3}\int   u^3\psi '-4\int   u^2\widehat{u}\psi '-\frac{1}{2}\int   u^2\widehat{u}\psi '''+\frac{1}{2}\int   u^2\widehat{u}_x\psi ''+\int   uh\psi '+\frac{1}{2}\int   uh_x\psi ''\\
&-\frac{5}{2}\int   \widehat{u}h_x\psi ''-2\int   \widehat{u}_xh\psi ''+\frac{1}{2}\int   \widehat{u}h_x\psi ^{(4)}+2\kappa \int \widehat{u}(1-\partial_x^2)^{-1}\widehat{u}_x(\psi ^{(4)}-5\psi '')\\
&-8\kappa\int \widehat{u}_x(1-\partial_x^2)^{-1}\widehat{u}_{xx}\psi ''-2\kappa \int y(1-\partial_x^2)^{-1}\widehat{u}_x\psi -2\kappa \int u_x\widehat{u}\psi \\
&-\dot{m_j}(t)\int (4\widehat{u}^2+5\widehat{u}_x^2+\widehat{u}_{xx}^2)\psi',
\end{split}
\end{align}
where $h\triangleq(1-\partial_x^2)^{-1}(u^2)$ and for notational convenience the argument $x-m_j(t)$ of $\psi$ is suppressed.
Consider first the quadratic terms, say the last fives. Denoting $r\triangleq(1-\partial_x^2)^{-1}\widehat{u},\quad \widehat{u}=r-r_{xx},$ and making use of the identity
\begin{align}\label{id1}
\begin{split}
\aaa \widehat{u}^2\psi'''=\aaa(r-r_{xx})^2\psi'''
=&\aaa (r^2+r_{xx}^2+2r_x^2)\psi'''-\aaa r^2\psi^{(5)},
\end{split}
\end{align}
we have
\begin{align}\label{bb3}
\begin{split}
\aaa \widehat{u}(1-\partial_x^2)^{-1}\widehat{u}_x(\psi^{(4)}-5\psi'')=&\aaa (r-r_{xx})r_x(\psi^{(4)}-5\psi'')\\
=&-\frac{1}{2}\aaa r^2(\psi^{(5)}-5\psi''')+\frac{1}{2}\aaa r_x^2(\psi^{(5)}-5\psi'''),
\end{split}
\end{align}
\begin{align}\label{bb2}
\begin{split}
\aaa \widehat{u}_x(1-\partial_x^2)^{-1}\widehat{u}_{xx}\psi''=&\aaa \widehat{u}_x(r-\widehat{u})\psi''=\aaa (r_x-r_{xxx})r\psi''+\frac{1}{2}\aaa \widehat{u}^2\psi'''\\
=&-\frac{1}{2}\aaa r^2\psi'''-\frac{3}{2}\aaa r_x^2\psi'''+\frac{1}{2}\aaa r^2\psi^{(5)}+\frac{1}{2}\aaa \widehat{u}^2\psi'''\\\overset{\eqref{id1}}{=}& \quad \frac{1}{2}\aaa (r_{xx}^2-r_x^2)\psi'''.
\end{split}
\end{align}
Integrating by parts leads to
\begin{align}\label{bb4}
\begin{split}
\aaa y(1-\partial_x^2)^{-1}\widehat{u}_x\psi=&\aaa [(1-\partial_x^2)(4-\partial_x^2)(1-\partial_x^2)r]r_x\psi\\
=&-2\aaa r^2\psi'+\frac{9}{2}\aaa r_x^2\psi'+9\aaa r_{xx}^2\psi'-3\aaa r_x^2\psi'''\\
&+\frac{5}{2}\aaa r_{xxx}^2\psi'-\frac{5}{2}\aaa r_{xx}^2\psi'''+\frac{1}{2}\aaa r_x^2\psi^{(5)},
\end{split}
\end{align}
\begin{align}\label{bb1}
\begin{split}
\aaa u_x\widehat{u}\psi=\aaa(4-\partial_x^2)\widehat{u}_x\widehat{u}\psi=&\aaa 4\widehat{u}_x\widehat{u}\psi-\aaa \widehat{u}_{xxx}\widehat{u} \psi\\
=&-2\aaa \widehat{u}^2\psi'-\frac{3}{2}\aaa \widehat{u}_x^2\psi'+\frac{1}{2}\aaa \widehat{u}^2\psi'''.
\end{split}
\end{align}
We compute $A\triangleq 2\kappa \cdot\eqref{bb3}-8\kappa\cdot\eqref{bb2}-2\kappa \cdot\eqref{bb4}-2\kappa \cdot\eqref{bb1}$:
\begin{align*}
\begin{split}
A
=&\kappa\aaa (4\widehat{u}^2+3\widehat{u}_x^2)\psi'-\kappa\aaa \widehat{u}^2\psi'''-\kappa\aaa r^2\psi^{(5)}\\
&+\kappa\aaa (4r^2-9r_x^2-18r_{xx}^2-5r_{xxx}^2)\psi'+\kappa\aaa (5r^2+5r_x^2+r_{xx}^2)\psi'''\\
=&\kappa\aaa (4\widehat{u}^2+3\widehat{u}_x^2)\psi'+\kappa\aaa (4r^2-9r_x^2-18r_{xx}^2-5r_{xxx}^2)\psi'+\kappa\aaa (4r^2+3r_x^2)\psi'''\quad\mbox{ By }\eqref{id1}\\
\leq&\kappa\aaa (4\widehat{u}^2+3\widehat{u}_x^2)\psi'+\kappa\aaa 4r^2(\psi'+\psi''')-\kappa\aaa r_x^2(\psi'-3\psi''')\\
\leq&\kappa\aaa (4\widehat{u}^2+3\widehat{u}_x^2)\psi'+\kappa(1+\frac{1}{B^2})\aaa 4r^2\psi'-\kappa(1-\frac{3}{B^2})\aaa r_x^2\psi'\\
\leq&\kappa\aaa (4\widehat{u}^2+3\widehat{u}_x^2)\psi'+\kappa(1+\frac{1}{B^2})\aaa 4r^2\psi',\quad\quad\mbox{by choosing}\quad B>2.
\end{split}
\end{align*}
We have as in \eqref{id1} that
\begin{align}\label{id2}
\begin{split}
\aaa \widehat{u}^2\psi'
=&\aaa (r^2+r_{xx}^2+2r_x^2)\psi'-\aaa r^2\psi'''\geq\aaa r^2(\psi'-\psi''')\geq(1-\frac{1}{B^2})\aaa r^2\psi'.
\end{split}
\end{align}
Choosing $B$ large such that $\frac{2B^2}{B^2-1}<2+\frac{\sigma_0}{\kappa}$, it follows that
\begin{align}
\begin{split}
A
\leq &\kappa\aaa (4\widehat{u}^2+3\widehat{u}_x^2)\psi'+\kappa(1+\frac{1}{B^2})(1-\frac{1}{B^2})^{-1}\aaa 4\widehat{u}^2\psi'=\frac{2B^2}{B^2-1}\kappa\aaa 4\widehat{u}^2\psi'+\kappa\aaa 3\widehat{u}_x^2\psi'\\
<&(2\kappa +\sigma_0)\aaa (4\widehat{u}^2+5\widehat{u}_x^2+\widehat{u}_{xx}^2)\psi'.
\end{split}
\end{align}
Then by \eqref{virialpsi} and the above, we deduce that
\begin{align}\label{326}
\begin{split}
\frac{d}{dt}\int (4\widehat{u}^2+5\widehat{u}_x^2+\widehat{u}_{xx}^2)\psi
\leq\ &\frac{2}{3}\int   u^3\psi '-4\int   u^2\widehat{u}\psi '-\frac{1}{2}\int   u^2\widehat{u}\psi '''+\frac{1}{2}\int   u^2\widehat{u}_x\psi ''+\int   uh\psi '\\
&+\frac{1}{2}\int   uh_x\psi ''-\frac{5}{2}\int   \widehat{u}h_x\psi ''-2\int   \widehat{u}_xh\psi ''+\frac{1}{2}\int   \widehat{u}h_x\psi ^{(4)}\\
&+(2\kappa +\sigma_0-\dot{m_j}(t))\int (4\widehat{u}^2+5\widehat{u}_x^2+\widehat{u}_{xx}^2)\psi\\
\triangleq&\sum_{k=1}^9I_j^k+(2\kappa +\sigma_0-\dot{m_j}(t))\int (4\widehat{u}^2+5\widehat{u}_x^2+\widehat{u}_{xx}^2)\psi\\
\leq&\sum_{k=1}^9I_j^k-\frac{\sigma_0}{2} \int (4\widehat{u}^2+5\widehat{u}_x^2+\widehat{u}_{xx}^2)\psi.
\end{split}
\end{align}
where \eqref{c1-0} and \eqref{c2-0} lead that $\dot{m_j}(t)>2\ka+\frac{3\sigma_0}{2}$. We claim that for $k=1,\cdots,9,$ it holds
\begin{equation}\label{JK}
I_j^k\leq \frac{\sigma_0}{20}\int (4\widehat{u}^2+5\widehat{u}_x^2+\widehat{u}_{xx}^2)\psi+Ce^{-\frac{1}{B}(\sigma_0t/2+L/4)}.
\end{equation}
Divide $\R$ into two regions $D_j$ and $D_j^c$ with $
D_j=\left[x_{j-1}(t)+\frac{L}{4},\  x_j(t)-\frac{L}{4}\right],\quad j=2,\cdots,N. $
Combining \eqref{c1-0} and \eqref{c2-0}, one check that for $x\in D_j^c$,
\begin{equation}
\vert x-m_j(t)\vert\geq\frac{x_j(t)-x_{j-1}(t)}{2}-\frac{L}{4}\geq \frac{\sigma_0t}{2}+\frac{x_j(0)-x_{j-1}(0)}{2}-\frac{L}{4}\geq \frac{\sigma_0t}{2}+\frac{L}{4}.
\end{equation}
Note that by Proposition \ref{hlinfty} and \eqref{exp0}:
\begin{align}\label{ubound}
\begin{split}
\Vert u\Vert_{L^\infty(D_j)}\leq&\left\Vert u-\sum\phi_{c_j(t)}(\cdot-x_j(t))\right\Vert_{L^\infty}+\sum\left\Vert \phi_{c_j(t)}(\cdot-x_j(t))\right\Vert_{L^\infty(D_j)}\\
=&O\left(\left\Vert u-\sum\phi_{c_j(t)}(\cdot-x_j(t))\right\Vert^{\frac{2}{3}}_{L^2}\right)+O\left(e^{-\sigma_0L/4}\right)=O(\al^{\frac{2}{3}})+O\left(e^{-\sigma_0L/4}\right),\\
\Vert u\Vert_{L^\infty}\quad\leq&\left\Vert u-\sum\phi_{c_j(t)}(\cdot-x_j(t))\right\Vert_{L^\infty}+\sum\left\Vert \phi_{c_j(t)}(\cdot-x_j(t))\right\Vert_{L^\infty}=O(\al^{\frac{2}{3}})+O(1),\\
\Vert \widehat{u}\Vert_{L^\infty(D_j)}
\leq&\left\Vert(4-\partial_x^2)^{-1}\left[ u-\sum\phi_{c_j(t)}(\cdot-x_j(t))\right]\right\Vert_{L^\infty}+\sum\left\Vert (4-\partial_x^2)^{-1}\left[\phi_{c_j(t)}(\cdot-x_j(t))\right]\right\Vert_{L^\infty(D_j)}\\
\leq&\left\Vert   u-\sum\phi_{c_j(t)}(\cdot-x_j(t)) \right\Vert_{L^2}+O(e^{-\sigma_0L/4})=O(\al)+O\left(e^{-\sigma_0L/4}\right).
\end{split}
\end{align}
One also has,
\begin{align}\label{l2bound}
\begin{split}
\aaa \vert u\vert(\vert\widehat{u}\vert+\vert\widehat{u}_x\vert+\vert\widehat{u}_{xx}\vert)\psi'=&\aaa \vert 4\widehat{u}-\widehat{u}_{xx}\vert (\vert\widehat{u}\vert+\vert\widehat{u}_x\vert+\vert\widehat{u}_{xx}\vert)\psi' \leq20\aaa \left(4\widehat{u}^2+5\widehat{u}_x^2+\widehat{u}_{xx}^2\right)\psi',\\
\aaa \vert u\vert(\vert\widehat{u}\vert+\vert\widehat{u}_x\vert+\vert\widehat{u}_{xx}\vert)=&\aaa \vert 4\widehat{u}-\widehat{u}_{xx}\vert (\vert\widehat{u}\vert+\vert\widehat{u}_x\vert+\vert\widehat{u}_{xx}\vert)\leq20\aaa \left(4\widehat{u}^2+5\widehat{u}_x^2+\widehat{u}_{xx}^2\right)\leq 20\Vert u\Vert_{L^2}^2,
\end{split}
\end{align}
and
\begin{align}\label{hbound}
\begin{split}
\aaa h\psi'
=&\aaa u^2\left[(1-\partial_x^2)^{-1}\psi'\right]
\leq\aaa u^2(1-\frac{1}{B^2})^{-1}\psi'\leq5\aaa \left(4\widehat{u}^2+5\widehat{u}_x^2+\widehat{u}_{xx}^2\right)\psi',\\
\aaa \vert u\vert h
=&\aaa \vert u\vert \left[(1-\partial_x^2)^{-1}(u^2)\right]
\leq\left\Vert (1-\partial_x^2)^{-1}\vert u\vert\ \right\Vert_{L^\infty}\aaa u^2
\leq\Vert u\Vert_{L^2}^3\\
\aaa \vert \widehat{u}\vert h
=&\aaa \left[(1-\partial_x^2)^{-1}\vert \widehat{u}\vert\right] u^2
\leq\left\Vert (1-\partial_x^2)^{-1}\vert \widehat{u}\vert\right\Vert_{L^\infty}\aaa u^2\leq \Vert u\Vert^3_{L^2},
\end{split}
\end{align}
where the first $\leq$ holds since $(1-\partial_x^2)\psi'=\psi'-\frac{1}{B^2}\psi'''\geq(1-\frac{1}{B^2})\psi',$  the last inequality holds because
\begin{align*}
(1-\partial_x^2)^{-1}\vert \widehat{u}(x)\vert
\leq&(1-\partial_x^2)^{-1}(4-\partial_x^2)^{-1}\vert u\vert
=\frac{1}{3}(1-\partial_x^2)^{-1}\vert u\vert-\frac{1}{3}(4-\partial_x^2)^{-1}\vert u\vert\\
\leq&\frac{1}{3}\left\Vert (1-\partial_x^2)^{-1}\vert u\vert\right\Vert_{L^\infty}+\frac{1}{3}\left\Vert (4-\partial_x^2)^{-1}\vert u\vert\right\Vert_{L^\infty} \quad \leq \quad \Vert u\Vert_{L^2}.
\end{align*}
One then has
\begin{align*}
I_j^1=&\frac{2}{3}\int_{D_j}u^3\psi'+\frac{2}{3}\int_{D_j^c}u^3\psi'
\leq\frac{2}{3}\Vert u\Vert_{L^\infty(D_j)}\aaa u^2\psi'+\frac{2}{3}\Vert\psi'\Vert_{L^\infty(D_j^c)}\Vert u\Vert_{L^3}\\
\leq&\frac{2}{3}\Vert u\Vert_{L^\infty(D_j)}\aaa \vert u(4\widehat{u}-\widehat{u}_{xx})\vert\psi'+Ce^{-\frac{1}{B}(\sigma_0t/2+L/4)}\Vert u\Vert_{L^\infty}\cdot\Vert u\Vert^2_{L^2},\\
I_j^2
\leq&4\Vert u\Vert_{L^\infty(D_j)}\aaa \vert u\widehat{u} \vert \psi'+Ce^{-\frac{1}{B}(\sigma_0t/2+L/4)}\Vert u\Vert_{L^\infty}\aaa \vert u\widehat{u}\vert, \\
I_j^5
\leq&\Vert u\Vert_{L^\infty(D_j)}\aaa  h\psi'+Ce^{-\frac{1}{B}(\sigma_0t/2+L/4)} \aaa \vert u\vert h,\\
I_j^7
\leq&\frac{5}{2}\Vert \widehat{u}\Vert_{L^\infty(D_j)}\aaa h\psi'+Ce^{-\frac{1}{B}(\sigma_0t/2+L/4)} \aaa \vert \widehat{u}\vert h.
\end{align*}
It follows from \eqref{ubound}, \eqref{l2bound}, \eqref{hbound} that for $\al$ small enough and $L$ large enough, $$I_j^{1,2,5,7}\leq \frac{\sigma_0}{20}\int (4\widehat{u}^2+5\widehat{u}_x^2+\widehat{u}_{xx}^2)\psi+Ce^{-\frac{1}{B}(\sigma_0t/2+L/4)}.$$
By \eqref{psi}, \eqref{ubound}, \eqref{l2bound}, the estimate for $I_j^{3,4}$ are almost the same as that for $I_j^2$.
Noticing that
\begin{align*}
h(x)=&\frac{e^{-x}}{2}\int_{-\infty}^x e^ru^2(r)dr+\frac{e^x}{2}\int_x^{\infty}e^{-r}u^2(r)dr,\\
h_x(x)=&\frac{-e^{-x}}{2}\int_{-\infty}^x e^ru^2(r)dr+\frac{e^x}{2}\int_x^{\infty}e^{-r}u^2(r)dr,
\end{align*}
we infer that $
\vert h_x(x)\vert\leq h(x),\quad \forall x\in\R.$
Then the estimate for $I_j^6$ follows from the same arguments as that of $I_j^5$. Noting that $I_j^8=2\aaa \widehat{u}h_x\psi''+2\aaa \widehat{u}h\psi'''$,
where the first integral is $I_j^7$, and the second could be  estimated as for $I_j^7$ using the inequality $\vert \psi'''\vert\leq\psi'$. The estimate for $I_j^9$ is exactly the same as for $I_j^7$ using $\vert\psi^{(4)}\vert\leq 3\psi'$.

It follows from \eqref{326} and \eqref{JK} that $
\frac{d}{dt}\int (4\widehat{u}^2+5\widehat{u}_x^2+\widehat{u}_{xx}^2)\psi \leq C e^{-\frac{1}{B}(\sigma_0t/2+L/4)}.$ The lemma follows from integrating the above from $0$ to $t$.
\end{proof}

\begin{proof}[Proof of Lemma \ref{Posi}]
Let $\eta>0$ be some small constant to be chosen later. Let $\Phi\in C^2(\R)$ be an even function, with
\begin{equation}
\Phi(x)=1-\eta_1\ \mbox{on}\ [0,1]; \ \Phi(x)=e^{-x}\ \mbox{on} \ [2,\infty);\  e^{-x}<\Phi(x)<3e^{-x} \ \mbox{on}\ (1,2); \ \mbox{and}\  \vert\Phi'\vert\leq\Phi.
\end{equation}
Let $\Phi_B(x)=\Phi(\frac{x}{B}).$
One has
\begin{equation}
\vert\Phi'_B(x)\vert=\frac{1}{B}\vert\Phi'(\frac{x}{B})\vert\leq\frac{1}{B}\Phi_B(x).
\end{equation}
Denoting $\Phi_{B_j}=\Phi_{B_j}(x,t)=\Phi_B(x-x_j(t)),$  it  thus transpires that
\begin{eqnarray*}
&&\int \left[-R(t)\e^2-2\kappa \e\widehat{\e}+c(t,x)(4\widehat{\e}^2+5\widehat{\e}_x^2+\widehat{\e}_{xx}^2)\right]dx\\
&=&\int \left[-R(t)\e^2-2\kappa (4\widehat{\e}^2+\widehat{\e}_x^2)+c(t,x)(4\widehat{\e}^2+5\widehat{\e}_x^2+\widehat{\e}_{xx}^2)\right]dx\\
&=&\sum \int \Phi_{B_j}\left[-R_j(t)\e(t)^2-2\kappa (4\widehat{\e}^2+\widehat{\e}_x^2)+c_j(t)(4\widehat{\e}^2+5\widehat{\e}_x^2+\widehat{\e}_{xx}^2)\right]dx\\
&&-\int \left(R(t)-\sum \Phi_{B_j}R_j(t)\right)\e^2(t) \ dx+\aaa\left(1-\sum \Phi_{B_j}\right)(c(t,x)-2\kappa )(4\widehat{\e}^2+\widehat{\e}_x^2)dx\\
&&+\aaa\left(1-\sum \Phi_{B_j}\right) c(t,x) (4\widehat{\e}_x^2+\widehat{\e}_{xx}^2)dx+\aaa \sum \Phi_{B_j}\left(c(t,x)-c_j(t)\right) (4\widehat{\e}^2+5\widehat{\e}_x^2+\widehat{\e}_{xx}^2)dx,
\end{eqnarray*}
where we applied the spectral information of $\cal$ in the last step.

Note that for $B$ sufficiently large, we have $\frac{1}{B}<\sigma_0$, which implies that the exponential decay rate of $\Phi_B$ is smaller than those of $\phi_{c_j(t)}$. Then from the orthogonality \eqref{ortho} and the exponential decay of $R_j$, it is easy to verify that for $B$ sufficiently large,
\begin{equation}
\left\vert\left(\e\sqrt{\Phi_{B_j}},R_j\right)_{\cas}\right\vert+\left\vert\left(\e\sqrt{\Phi_{B_j}},R_{j,x}\right)_{\cas}\right\vert<\eta\left\Vert \e\sqrt{\Phi_{B_j}}\right\Vert_{L^2}.
\end{equation}
Then by Lemma \ref{Positivity}, we have
\begin{align}\label{123}
\begin{split}
\int -R_j(t)&\e(t)^2\Phi_{B_j}-2\kappa \e\sqrt{\Phi_{B_j}}(4-\partial_x^2)^{-1}\left(\e\sqrt{\Phi_{B_j}}\right) \\
&\quad+c_j(t) \left(\e\sqrt{\Phi_{B_j}}\right)(1-\partial_x^2)(4-\partial_x^2)^{-1}\left(\e\sqrt{\Phi_{B_j}}\right)  dx
 \geq\ \theta\int \e^2\Phi_{B_j}
\end{split}
\end{align}
One has
\begin{align}
\begin{split}
&\left\vert\aaa\left[\e\sqrt{\Phi_{B_j}(x)}(4-\partial_x^2)^{-1}\left(\e\sqrt{\Phi_{B_j}}\right)(x)-\left(\e \Phi_{B_j}\right)(x)(4-\partial_x^2)^{-1} \e
\right]dx\right\vert\\
=&\frac{1}{4}\left\vert\aaa \e(x)\aaa e^{-2\vert y-x\vert}\e(y)\sqrt{\Phi_{B_j}(x)}\left(\sqrt{\Phi_{B_j}(x)}-\sqrt{\Phi_{B_j}(y)}\right)dydx\right\vert\\
\leq&\frac{\Vert\e\Vert_{L^2}}{4}\sqrt{\aaa\left(\aaa e^{-2\vert y-x\vert}\e(y)\sqrt{\Phi_{B_j}(x)}\left(\sqrt{\Phi_{B_j}(x)}-\sqrt{\Phi_{B_j}(y)}\right)dy\right)^2dx}\\
\leq&\frac{\Vert\e\Vert^2_{L^2}}{4}\sqrt{\aaa\aaa e^{-4\vert y-x\vert} \Phi_{B_j}(x)\left(\sqrt{\Phi_{B_j}(x)}-\sqrt{\Phi_{B_j}(y)}\right)^2dydx}\\
=&\frac{\Vert\e\Vert^2_{L^2}}{4}\sqrt{\aaa\aaa e^{-4\vert y-x\vert} \Phi_B(x)\left(\sqrt{\Phi_B(x)}-\sqrt{\Phi_B(y)}\right)^2dydx},
\end{split}
\end{align}
and
\begin{align}
\begin{split}
\aaa\aaa e^{-4\vert y-x\vert} \Phi_B(x)&\left(\sqrt{\Phi_B(x)}-\sqrt{\Phi_B(y)}\right)^2dydx\\
=& \aaa\aaa e^{-4\vert y-x\vert} \Phi_B(x)\left(\frac{\Phi'_B}{2\sqrt{\Phi_B}}\cdot(x-y)\right)^2dydx\\
\leq& \aaa\aaa e^{-4\vert y-x\vert} \Phi_B(x)\left(\frac{\Phi_B}{2B\sqrt{\Phi_B}}(x-y)\right)^2dydx\\
=&\frac{1}{4B^2}\aaa\aaa e^{-4\vert y-x\vert} \Phi_B(x) (x-y) ^2dydx\quad
=\quad\frac{O(1)}{B},
\end{split}
\end{align}
which implies
\begin{equation}
\left\vert\aaa\left(\e\sqrt{\Phi_{B_j}(x)}(4-\partial_x^2)^{-1}\left(\e\sqrt{\Phi_{B_j}}\right)(x)-\left(\e \Phi_{B_j}\right)(x)(4-\partial_x^2)^{-1} \e
\right)dx\right\vert\leq\frac{O(1)}{\sqrt{B}}\Vert\e\Vert_{L^2}^2.
\end{equation}
It follows that
\begin{align}\label{124}
\begin{split}
\int -R_j(t)&\e(t)^2\Phi_{B_j}-2\kappa \e\sqrt{\Phi_{B_j}}(4-\partial_x^2)^{-1}\left(\e\sqrt{\Phi_{B_j}}\right) \\
&\quad+c_j(t) \left(\e\sqrt{\Phi_{B_j}}\right)(1-\partial_x^2)(4-\partial_x^2)^{-1}\left(\e\sqrt{\Phi_{B_j}}\right)  dx \\
=&\int (c_j(t)-R_j(t))\e(t)^2\Phi_{B_j}-(3c_j(t)+2\kappa )\e\sqrt{\Phi_{B_j}}(4-\partial_x^2)^{-1}\left(\e\sqrt{\Phi_{B_j}}\right) \\
=&\int (c_j(t)-R_j(t))\e(t)^2\Phi_{B_j}-(3c_j(t)+2\kappa )\e\Phi_{B_j}(4-\partial_x^2)^{-1} \e  +\frac{O(1)}{\sqrt{B}}\Vert\e\Vert_{L^2}^2 \\
=&\int -R_j(t)\e(t)^2\Phi_{B_j}-2\kappa \Phi_{B_j}(4\widehat{\e}^2+\widehat{\e}_x^2)+c_j\Phi_{B_j}(4\widehat{\e}^2+5\widehat{\e}_x^2+\widehat{\e}_{xx}^2) +\frac{O(1)}{\sqrt{B}}\Vert\e\Vert_{L^2}^2, \\
\end{split}
\end{align}
where in the last step, terms containing $\Phi_{B_j}'$ from integration by parts  are of  $\frac{O(1)}{B}\Vert\e\Vert_{L^2}^2$, which are controlled by $\frac{O(1)}{\sqrt{B}}\Vert\e\Vert_{L^2}^2$.
\eqref{124} and \eqref{123} imply
\begin{align}
\begin{split}
\sum &\int \Phi_{B_j}\left[-R_j(t)\e(t)^2 -2\kappa (4\widehat{\e}^2+\widehat{\e}_x^2)+c_j(t)(4\widehat{\e}^2+5\widehat{\e}_x^2+\widehat{\e}_{xx}^2)\right]dx\\
&\geq\theta\sum \int \Phi_{B_j} \e(t)^2 dx-\frac{O(1)}{\sqrt{B}}\Vert\e\Vert_{L^2}^2
\geq\theta\sum \int \Phi_{B_j} (4\widehat{\e}^2+5\widehat{\e}_x^2+\widehat{\e}_{xx}^2) dx-\frac{O(1)}{\sqrt{B}}\Vert\e\Vert_{L^2}^2.
\end{split}
\end{align}
Since $1-\sum \Phi_{B_j}>0,\quad c(t,x)-2\kappa >\sigma_0, \quad c(t,x)>\sigma_0,$
\begin{align}
\begin{split}
\aaa\left(1-\sum \Phi_{B_j}\right)(c(t,x)-2\kappa )&(4\widehat{\e}^2+\widehat{\e}_x^2)dx+\aaa\left(1-\sum \Phi_{B_j}\right) c(t,x) (4\widehat{\e}_x^2+\widehat{\e}_{xx}^2)dx\\
&\geq \aaa\left(1-\sum \Phi_{B_j}\right)\sigma_0(4\widehat{\e}^2+5\widehat{\e}_x^2+\widehat{\e}_{xx}^2)dx
\end{split}
\end{align}
Notice that $ \displaystyle \eta_1R(t)-\sum \Psi_j\phi_{c_j(t)}(x-x_j(t))$ is all exponentially small, being $O\left(e^{-\frac{1}{B}(\frac{\sigma_0t}{2}+\frac{L}{4})}+e^{-\sigma_0L/4}\right)$. If we choose $L_4=MB$,  with $M\gg 0$, and $\lambda_0=\lambda_1/5$, it then follows that
\begin{align}
\begin{split}
&\left\vert-\int \left(R(t)-\sum \Phi_{B_j}R_j(t)\right)\e^2(t) \ dx\right\vert\\
\leq &\eta_1\int  R(t) \e^2(t) \ dx+\int \left((1-\eta_1)R(t)-\sum \Phi_{B_j}R_j(t)\right)\e^2(t) \ dx=\eta_1O(1)\Vert \e\Vert_{L^2}^2.
\end{split}
\end{align}
Denote $\Psi_1(t,x)=1-\psi(x-m_2(t))$, $\Psi_j(t,x)=\psi(x-m_{j}(t))-\psi(x-m_{j+1}(t)),\ j=2,\cdots,N-1,$ and $\Psi_N(t,x)=\psi(x-m_N(t)).$ Then  $c(t,x)=\sum \Psi_j(t,x)c_j(t)$ and $\sum \Psi_j(t,x)=1.$
From the property of $\psi$, the fact that $\Phi_{B_j}\leq 3e^{-\frac{\vert x-m_j(t)\vert}{B}}$, and that $\vert m_j(t)-x_j(t)\vert\geq L_4/2\geq 2MB$, we obtain
\begin{align}
\begin{split}
\vert \Phi_{B_j}\left(c(t,x)-c_j(t)\right)\vert&\leq\vert c(t,x)-c_j(t)\vert_{L^\infty(\vert x-x_j(t)\vert\leq MB)}+Ce^{-M}\leq Ce^{-\sigma_0MB}+Ce^{-M}.
\end{split}
\end{align}
It follows from that for sufficiently small $\eta_1>0$,
\begin{align}
\begin{split}
&\int \left[-R(t)\e^2-2\kappa \e\widehat{\e}+c(t,x)(4\widehat{\e}^2+5\widehat{\e}_x^2+\widehat{\e}_{xx}^2)\right]dx\\
\geq&\theta\sum \int \Phi_{B_j} (4\widehat{\e}^2+5\widehat{\e}_x^2+\widehat{\e}_{xx}^2) dx+\aaa\left(1-\sum \Phi_{B_j}\right)\sigma_0(4\widehat{\e}^2+5\widehat{\e}_x^2+\widehat{\e}_{xx}^2)dx\\
&-\frac{O(1)}{\sqrt{B}}\Vert\e\Vert_{L^2}^2-\eta_1O(1)\Vert \e\Vert_{L^2}^2-(Ce^{-\sigma_0MB}+Ce^{-M})\aaa  (4\widehat{\e}^2+5\widehat{\e}_x^2+\widehat{\e}_{xx}^2)dx\\
\geq&\frac{\min\{\theta,\sigma_0\}}{2}\aaa  (4\widehat{\e}^2+5\widehat{\e}_x^2+\widehat{\e}_{xx}^2)dx-\frac{\min\{\theta,\sigma_0\}}{16}\Vert \e\Vert_{L^2}^2\quad
\geq\quad\frac{\min\{\theta,\sigma_0\}}{16}\Vert \e\Vert_{L^2}^2
\end{split}
\end{align}
\end{proof}

\begin{proof}[Proof of Lemma \ref{Qua}]
Noting that
\begin{equation*}
H(\phi_{c_j(t)})-H(\ppo)=(\frac{dH(\ppo)}{dc})(c_j(t)-c_j(0))+O(\vert c_j(t)-c_j(0)\vert^2),
\end{equation*}
\begin{equation*}
S(\phi_{c_j(t)})-S(\ppo)=(\frac{dS(\ppo)}{dc})(c_j(t)-c_j(0))+O(\vert c_j(t)-c_j(0)\vert^2),
\end{equation*}
it follows from \eqref{derivative} that
\begin{equation}\label{HS}
H(\phi_{c_j(t)})-H(\ppo)=-c_j(0)\left(
S(\phi_{c_j(t)})-S(\ppo)\right)+O(\vert c_j(t)-c_j(0)\vert^2).
\end{equation}

\noindent {\it Step 1.}  We prove that there exists $C>0$ such that
\begin{align}\label{42}
\begin{split}
\left\vert\sum_{j=1}^Nc_j(0)[S(\phi_{c_j(t)})-S(\phi_{c_j(0)})]\right\vert\leq &C\left(\Vert\e(t)\Vert_{L^2}^2 +\Vert\e(0)\Vert_{L^2}^2+\Vert\e(0)\Vert_{L^3}^3+e^{-\frac{\sigma_0L}{2}}+\sum_{j=1}^N\vert c_j(t)-c_j(0)\vert^2\right).
\end{split}
\end{align}
By \eqref{16} and \eqref{co2}, it is found that
\begin{align*}
\left\vert\sum_{j=1}^NH(\phi_{c_j(t)})-H(\phi_{c_j(0)})\right\vert
\leq& C\left(\Vert\e(t)\Vert_{L^3}^3+\Vert\e(t)\Vert_{L^2}^2+\Vert\e(0)\Vert_{L^2}^2+\Vert\e(0)\Vert_{L^3}^3+e^{-\frac{\sigma_0L}{2}}\right)\\
\leq& C\left((\Vert\e(t)\Vert_{L^\infty}+1)\Vert\e(t)\Vert_{L^2}^2+\Vert\e(0)\Vert_{L^2}^2+\Vert\e(0)\Vert_{L^3}^3+e^{-\frac{\sigma_0L}{2}}\right)
\end{align*}
This together with \eqref{HS} and \eqref{eeee} implies \eqref{42}.

\noindent {\it Step 2.}  We prove that for $
d_j(t)\triangleq\sum_{k=j}^NS(\phi_{c_j(t)})$,
\begin{equation}\label{j2}
d_j(t)-d_j(0)\leq-(d_j(t)-d_j(0))+C(\Vert\e(0)\Vert_{L^2}^2+e^{-\frac{1}{B}(\frac{\sigma_0t}{2}+\frac{L}{8})}+e^{-\sigma_0L/8}), \quad j\geq 2;
\end{equation}
\begin{equation}\label{j1}
\vert d_1(t)-d_1(0)\vert\leq\frac{1}{2}\Vert\e(0)\Vert_{\cas}^2-\frac{1}{2}\Vert\e(t)\Vert_{\cas}^2+Ce^{-\sigma_0L/8}.
\end{equation}
In order to prove \eqref{j2},  we shall prove
\begin{equation}\label{ddd}
\left\vert
I_j(t)-d_j(t)-\frac{1}{2}\int \psi(\cdot-m_j(t))(4\widehat{\e}^2(t)+5\widehat{\e}^2_x(t)+\widehat{\e}_{xx}^2(t))dx
\right\vert\leq C(e^{-\frac{1}{B}(\frac{\sigma_0t}{2}+\frac{L}{8})}+e^{-\sigma_0L/8}).
\end{equation}
Indeed, $
d_j(t)=\frac{1}{2}\sum_{l=j}^N\int\left(4\widehat{R}_l^2+5\left(\widehat{R}_{l,x}\right)^2+\left(\widehat{R}_{l,xx}\right)^2\right)dx$, and
\begin{align*}
&
2I_j(t)-2d_j(t)-\int \psi(\cdot-m_j(t))(4\widehat{\e}^2(t)+5\widehat{\e}^2_x(t)+\widehat{\e}_{xx}^2(t))dx\\
=&
\int \psi(\cdot-m_j(t))\left((4\widehat{u}^2+5\widehat{u}_x^2+\widehat{u}_{xx}^2)-(4\widehat{\e}^2(t)+5\widehat{\e}^2_x(t)+\widehat{\e}_{xx}^2(t))\right)-\sum_{l=j}^N\int \left(4\widehat{R}_l^2+5\widehat{R}_{l,x}^2+\widehat{R}_{l,xx}^2dx\right)\\
=&
\int \psi(\cdot-m_j(t))\left(\sum_{l=1}^N\left(4\widehat{R}_l^2+5\widehat{R}_{l,x}^2+\widehat{R}_{l,xx}^2+2(4\widehat{\e}\widehat{R_l}+5\widehat{\e}_x\widehat{R}_{l,x}+\widehat{\e}_{xx}\widehat{R}_{l,xx})\right)+\mbox{mixed terms}\right)dx\\
&
-\sum_{l=j}^N\int \left(4\widehat{R}_l^2+5\widehat{R}_{l,x}^2+\widehat{R}_{l,xx}^2dx\right)\\
=&
\int \psi(\cdot-m_j(t))\left(\sum_{l=1}^{j-1}\left(4\widehat{R}_l^2+5\widehat{R}_{l,x}^2+\widehat{R}_{l,xx}^2+2(4\widehat{\e}\widehat{R_l}+5\widehat{\e}_x\widehat{R}_{l,x}+\widehat{\e}_{xx}\widehat{R}_{l,xx})\right)\right)dx\\
&
-\int \left(1-\psi(\cdot-m_j(t))\right)\sum_{l=j}^{N}\left(4\widehat{R}_l^2+5\widehat{R}_{l,x}^2+\widehat{R}_{l,xx}^2+2(4\widehat{\e}\widehat{R_l}+5\widehat{\e}_x\widehat{R}_{l,x}+\widehat{\e}_{xx}\widehat{R}_{l,xx})\right)+O\left(e^{-\sigma_0L/8}\right),\\
\end{align*}
where we have used the first orthogonality of \eqref{ortho}, which is equivalently expressed as
\begin{equation}\label{1111}
\aaa \left(4\widehat{\e}\widehat{R_l}+5\widehat{\e}_x\widehat{R}_{l,x}+\widehat{\e}_{xx}\widehat{R}_{l,xx}\right)dx=0.
\end{equation}
To evaluate the first integral, we divide $\R$ into two intervals by the point $x_{j-1}(t)+L/8$. To the left of the point,  \begin{eqnarray*}
\vert x-m_j(t)\vert
&\geq&\frac{x_j(t)-x_{j-1}(t)}{2}-\frac{L}{8}\geq\frac{\sigma_0t}{2}+\frac{x_j(0)-x_{j-1}(0)}{2}-\frac{L}{8}\geq\frac{\sigma_0t}{2}
+\frac{L}{8},\\
\vert\psi(x-m_j(t))\vert
&\leq& Ce^{-\frac{1}{B}(\frac{\sigma_0t}{2}+\frac{L}{8})},\quad\quad\quad\quad\quad\quad\quad
\end{eqnarray*}
and to the right of the point, $R_l$, $l\leq j-1$ are exponentially small, being $O\left(e^{-\sigma_0L/8}\right)$.
To evaluate the second integral, we divide $\R$ into two intervals by the point $x_{j}(t)-L/8$. To the right of the point,  \begin{eqnarray*}
\vert x-m_j(t)\vert
&\geq&\frac{x_j(t)-x_{j-1}(t)}{2}-\frac{L}{8}\geq\frac{\sigma_0t}{2}+\frac{x_j(0)-x_{j-1}(0)}{2}-\frac{L}{8}\geq\frac{\sigma_0t}{2}
+\frac{L}{8},\\
\vert1-\psi(x-m_j(t))\vert
&\leq& Ce^{-\frac{1}{B}(\frac{\sigma_0t}{2}+\frac{L}{8})},
\end{eqnarray*}
and to the left of the point, $R_l$, $l\geq j$ are exponentially small, being $O\left(e^{-\sigma_0L/8}\right)$. This proves \eqref{ddd}.

Then Lemma \ref{monotone} together with \eqref{ddd} imply
\begin{equation}\label{e:19}
\begin{aligned}
d_j(t)-d_j(0)\leq&\frac{1}{2}\int_{\R}\psi(\cdot-m_j(s))(4\hat{\e}^2(s)+5\hat{\e}^2_x(s)+\hat{\e}_{xx}^2(s))dx\bigg\vert_{s=t}^{s=0}
+C\left(e^{-\frac{1}{B}(\frac{\sigma_0t}{2}+\frac{L}{8})}+e^{-\sigma_0L/8}\right).
\end{aligned}
\end{equation}
\eqref{j2} then follows from the above estimate easily.

To prove \eqref{j1}, we note that $S(u(t))=S(u(0))$ and by orthogonality \eqref{1111} that
\begin{align}
\begin{split}
\frac{1}{2}\aaa (4\widehat{u}^2+5\widehat{u}_x^2+\widehat{u}_{xx}^2)=&\frac{1}{2}\aaa (4\widehat{R}^2+5\widehat{R}_x^2+\widehat{R}_{xx}^2)+\frac{1}{2}\aaa (4\widehat{\e}^2+5\widehat{\e}_x^2+\widehat{\e}_{xx}^2)\\
=&d_1(t)+O(e^{-\frac{\sigma_0L}{2}})+ S(\e(t)).
\end{split}
\end{align}
Then \eqref{j1} follows.

{\it Step 3. } Resummation by the Abel transform:
\begin{align*}
\sum_{j=1}^Nc_j(0)&\left[S(\phi_{c_j(t)})-S(\phi_{c_j(0)})\right] \\
=& \sum_{j=1}^{N-1}c_j(0)[(d_j(t)-d_{j+1}(t))-(d_j(0)-d_{j+1}(0))] +c_N(0)[d_N(t)-d_N(0)]\\
=&c_1(0)[d_1(t)-d_1(0)]+\sum_{j=2}^{N}[c_j(0)-c_{j-1}(0)][d_j(t)-d_j(0)].
\end{align*}
It then follows from  \eqref{42}  that
\begin{align}\label{50}
\begin{split}
-c_1(0)[d_1(t)-d_1(0)]-&\sum_{j=2}^{N}[c_j(0)-c_{j-1}(0)][d_j(t)-d_j(0)]\\
&\leq C\left(\Vert\e(t)\Vert_{L^2}^2+\Vert\e(0)\Vert_{L^2}^2+e^{-\sigma_0L/4}+\sum_{j=1}^N\vert c_j(t)-c_j(0)\vert^2\right).
\end{split}
\end{align}
Since $c_1(0)\geq\sigma_0,\ c_j(0)-c_{j-1}(0)\geq\sigma_0$, by \eqref{j1} and \eqref{j2}, we have
\begin{eqnarray*}
\sigma_0\sum_{j=1}^N\vert d_j(t)-d_j(0)\vert
&\leq& c_1(0)\vert d_1(t)-d_1(0)\vert+\sum_{j=2}^N[c_j(0)-c_{j-1}(0)] \ \vert d_j(t)-d_j(0)\vert\\
&\leq&-\left(c_1(0)[d_1(t)-d_1(0)]+\sum_{j=2}^{N}[c_j(0)-c_{j-1}(0)][d_j(t)-d_j(0)]\right)\\
&&+C\left(\Vert\e(0)\Vert_{L^2}^2+e^{-\frac{1}{B}(\frac{\sigma_0t}{2}+\frac{L}{8})}+e^{-\sigma_0L/8}\right).
\end{eqnarray*}
Then by \eqref{50}, we have
$$
\sum_{j=1}^N\vert d_j(t)-d_j(0)\vert
\leq C\left(\Vert\e(t)\Vert_{L^2}^2 +\Vert\e(0)\Vert_{L^2}^2+e^{-\frac{1}{B}(\frac{\sigma_0t}{2}+\frac{L}{8})}+e^{-\sigma_0L/8}+\sum_{j=1}^N\vert c_j(t)-c_j(0)\vert^2\right).
$$
Note that from the boundedness (from both above and below \eqref{derivative}) of $\frac{dS(\phi_c)}{dc}$, we have
$$
\vert c_j(t)-c_j(0)\vert \leq  C\vert S(\ppt)-S(\ppo)\vert\leq  C\left(\vert d_j(t)-d_j(0)\vert+\vert d_{j+1}(t)-d_{j+1}(0)\vert\right),
$$
which then yields
$$
\sum_{j=1}^N\vert c_j(t)-c_j(0)\vert
\leq C\left(\Vert\e(t)\Vert_{L^2}^2 +\Vert\e(0)\Vert_{L^2}^2+e^{-\frac{1}{B}(\frac{\sigma_0t}{2}+\frac{L}{8})}+e^{-\sigma_0L/8}
+\sum_{j=1}^N\vert c_j(t)-c_j(0)\vert^2\right).
$$
By choosing $\al$ small and $L$ large enough, we assume $C\vert c_j(t)-c_j(0)\vert\leq 1/2$ and so \eqref{quadratic} follows and the proof of the lemma is complete.
\end{proof}
\begin{proof}[Proof of Lemma \ref{Lya}]
Recall that $\gamma_0=\min\{1/(8B), \sigma_0/8\}.$ In view of \eqref{HS}, \eqref{j2} and  \eqref{j1},  we infer that
\[
\begin{aligned}
&-\int_{\R}\left[\frac{1}{2}R(t)\e^2(t)+k\e(t)(4-\partial_x^2)^{-1}\e(t)\right] dx\\
\overset{\eqref{16}}{\leq}&
-\sum_{j=1}^N\left[H(\ppt)-H(\ppo)\right]+ K_2\left\{\Vert\e(t)\Vert_{L^2}^2\Vert\e(t)\Vert_{L^\infty}+\Vert\e(0)\Vert_{L^2}^2+e^{-\sigma_0L/8}\right\}\\
\overset{\qquad}{\leq}&
\sum_{j=1}^Nc_j(0)\left[S(\ppt)-S(\ppo)\right] +C\sum_{j=1}^N\left(c_j(t)-c_j(0)\right)^2+C\left(\Vert\e(t)\Vert_{L^2}^{\frac{8}{3}}+\Vert\e(0)\Vert_{L^2}^2 +e^{-\gamma_0L}\right)\\
\overset{\eqref{quadratic}}{\leq}&c_1(0)[d_1(t)-d_1(0)]+\sum_{j=2}^{N}[c_j(0)-c_{j-1}(0)][d_j(t)-d_j(0)]+C\left(\Vert\e(t)\Vert_{L^2}^{\frac{8}{3}}+\Vert\e(0)\Vert_{L^2}^2 +e^{-\gamma_0L}\right)\\
\overset{\eqref{e:19}}{\leq}&\overset{\eqref{j1}}{\quad -}
c_1(0)[\frac{1}{2}\Vert\e(t)\Vert_{\cas}^2-\frac{1}{2}\Vert\e(0)\Vert_{\cas}^2 ]+\\
&\sum_{j=2}^{N}\left[c_j(0)-c_{j-1}(0)\right]\frac{1}{2}\left[   \int_{\R}\psi(\cdot-m_j(0))(4\hat{\e}^2(0)+5\hat{\e}^2_x(0)+\hat{\e}_{xx}^2(0))dx-\right.\\
&\left.\int_{\R}\psi(\cdot-m_j(t))(4\hat{\e}^2(t)+5\hat{\e}^2_x(t)+\hat{\e}_{xx}^2(t))dx \right] +C\left(\Vert\e(t)\Vert_{L^2}^{\frac{8}{3}}+\Vert\e(0)\Vert_{L^2}^2 +e^{-\gamma_0L}\right)\\
\overset{\qquad}{\leq}&
-c_1(0)\frac{1}{2}\Vert\e(t)\Vert_{\cas}^2-\sum_{j=2}^{N}[c_j(0)-c_{j-1}(0)]\frac{1}{2}\int_{\R}\psi(\cdot-m_j(t))(4\hat{\e}^2(t)+5\hat{\e}^2_x(t)+\hat{\e}_{xx}^2(t))dx+\\
&C\left(\Vert\e(t)\Vert_{L^2}^{\frac{8}{3}}+\Vert\e(0)\Vert_{L^2}^2 +e^{-\gamma_0L}\right)\\
\overset{\eqref{quadratic}}{\leq}&
-c_1(t)\frac{1}{2}\Vert\e(t)\Vert_{\cas}^2-\sum_{j=2}^{N}[c_j(t)-c_{j-1}(t)]\frac{1}{2}\int_{\R}\psi(\cdot-m_j(t))(4\hat{\e}^2(t)+5\hat{\e}^2_x(t)+\hat{\e}_{xx}^2(t))dx+\\
&C\left(\Vert\e(t)\Vert_{L^2}^{\frac{8}{3}}+\Vert\e(0)\Vert_{L^2}^2 +e^{-\gamma_0L}\right)\\
\overset{\eqref{positivity}}{=}&-\int_{\R}\frac{1}{2}c(t,x)\left(4\hat{\e}^2(t)+5\hat{\e}^2_x(t)+\hat{\e}_{xx}^2(t)\right)\ dx+C\left(\Vert\e(t)\Vert_{L^2}^{\frac{8}{3}}+\Vert\e(0)\Vert_{L^2}^2 +e^{-\gamma_0L}\right).
\end{aligned}
\]
Therefore,
\begin{equation*}
\int-R(t)\e^2(t)-2\kappa \e(t)\widehat{\e}(t)+c(t,x)\left(4\widehat{\e}^2(t)+5\widehat{\e}^2_x(t)+\widehat{\e}_{xx}(t)\right)
 \leq C\left(\Vert\e(t)\Vert_{L^2}^{\frac{8}{3}}+\Vert\e(0)\Vert_{L^2}^2 +e^{-\gamma_0L}\right).
\end{equation*}
Then from Lemma \ref{Posi}, we have $\lambda_0\Vert\e(t)\Vert_{L^2}^2
 \leq C\left(\Vert\e(t)\Vert_{L^2}^{\frac{8}{3}}+\Vert\e(0)\Vert_{L^2}^2 +e^{-\gamma_0L}\right),$ which implies \eqref{lya}.
\end{proof}

\noindent{\bf Acknowledgments}  The work of Li is partially supported by the NSFC grant 11771161. The work of Liu is partially supported by the Simons Foundation  grant 499875. The work of Wu is partially supported by the NSF grant DMS-1815079.


\bibliographystyle{plain}

\end{document}